\documentclass{birkjour}
\allowdisplaybreaks[4]

\input{mathrsfs.sty}

\usepackage[english]{babel}

\usepackage{amsmath,amssymb,amsthm,mathtools,enumerate,hyperref}

 \theoremstyle{plain}
 \newtheorem{thm}{Theorem}[section]
 \newtheorem{prob}{Problem}[section]
 
 \theoremstyle{definition}
 \newtheorem{defn}{Definition}[section]

 \theoremstyle{remark}
 \newtheorem{rem}{Remark}[section]
 
 \numberwithin{equation}{section}

\DeclareMathOperator{\td}{d\!}
\DeclareMathOperator{\te}{e}
\DeclareMathOperator{\arctanh}{arctanh}
\DeclareMathOperator{\arccot}{arccot}
\newcommand{\E}{\mathscr{E}}
\DeclareMathOperator{\bell}{B}
\DeclareMathOperator{\ti}{i}

\newcounter{main}
\renewcommand{\thesection}{\themain.\arabic{section}}
\renewcommand{\thesubsection}{\thesection.\arabic{subsection}}
\renewcommand{\thesubsubsection}{\thesubsection.\arabic{subsubsection}}

\begin{document}

\title[Theory of Normalized Remainders]{Chapter~12\\[1em]
Theory of Normalized Remainders in Taylor Series Expansions}

\author[F. Qi]{Feng Qi}

\address{School of Mathematics and Physics\\
Hulunbuir University\\
Hailar 021008, Inner Mongolia\\
China\\[0.5em]
17709 Sabal Court\\
University Village\\
Dallas, TX 75252-8024\\
USA}
\email{\href{mailto: F. Qi <qifeng618@gmail.com>}{qifeng618@gmail.com}}
\urladdr{\url{https://orcid.org/0000-0001-6239-2968}}

\subjclass{Primary 41A80; Secondary 11B68, 11B73, 11B75, 11B83, 11C08, 41A58}

\keywords{normalized remainder, definition, example, property, origin, background, theory of normalized remainder, Taylor series expansion, Maclaurin series expansion, conjecture, problem}

\begin{abstract}
Since 2023, through the detailed examination of numerous concrete examples, the author and his collaborators have identified a recurring pattern. Building upon this observation, they introduced the concept of the normalized remainder. They deliberately chose this term and subsequently explored its historical background and mathematical significance.
\par
In 2026, Abu-Ghuwaleh propelled the subject forward at a deeper structural level. By exploring the broader dynamical and theoretical framework surrounding the normalized remainder family, he significantly developed and formalized the concept, firmly embedding it within the field. Consequently, the notion of the normalized remainder now carries richer and more profound mathematical significance.
\par
In this chapter, the author presents a synthesis of the research process and the principal findings related to the normalized remainder.
\end{abstract}

\thanks{The author was partially supported by the Youth Project of Hulunbuir City for Basic Research and Applied Basic Research (Grant No.~GH2024020) and by the Natural Science Foundation of Inner Mongolia Autonomous Region (Grant No.~2025QN01041).}

\thanks{Please cite this article as ``Feng Qi, \textit{Theory of normalized remainders in Taylor series expansions}, Chapter~12 in Contemporary Approaches to Approximation of Operators---Approximation Theory and Mathematical Methods, edited by Gradimir V. Milovanovi\'c and Vijay Gupta, published in Birkh\"auser's \href{https://link.springer.com/series/4961}{Trends in Mathematics} series, pp.~349\nobreakdash--400, November 2026. URL: \url{https://link.springer.com/book/9783032387660}.''}

\maketitle
\tableofcontents

\section{Definition of the Normalized Remainder}
In this chapter, we use $\mathbb{N}$ and $\mathbb{N}_0$ to denote the sets of all positive integers and all nonnegative integers, respectively.
\par
What is the normalized remainders? The concept of normalized remainders in mathematics—particularly in the context of power series expansions—is indeed quite specialized and not widely familiar, even among mathematicians, because it is an innovation that emerged in 2023. In the context of the Maclaurin or Taylor series expansions, a remainder refers to the error term left after truncating the series at a certain degree. A normalized remainder is an essentially modified version of this error term.

\begin{defn}\label{defn-normaliz}
Let $f$ be an infinitely differentiable real-valued function on an interval $I\subseteq\mathbb{R}$ with $x_0\in I$ being an inner point. The formal series expansion of $f$ around the point $x_0$ is
\begin{equation}\label{Formal-series}
\sum_{n=0}^\infty f^{(n)}(x_0)\frac{(x-x_0)^n}{n!}.
\end{equation}
Its partial sum of the first $n+1$ terms is
\begin{equation*}
  S_n(x,x_0)=\sum_{k=0}^{n} f^{(k)}(x_0)\frac{(x-x_0)^k}{k!}, \quad n\in\mathbb{N}_0.
\end{equation*}
\begin{enumerate}[1)]
\item
If $f(x_0)\ne0$, then we define
\begin{equation*}
T_{-1}[f(x),x_0]=\frac{f(x)}{f(x_0)}, \quad x\in I.
\end{equation*}
\item
If $f^{(n+1)}(x_0)\ne0$ for $n\in\mathbb{N}_0$, then we define
\begin{equation}\label{T(f(x))}
T_n[f(x),x_0]=
\begin{dcases}
\frac1{f^{(n+1)}(x_0)}\frac{(n+1)!}{(x-x_0)^{n+1}}[f(x)-S_n(x,x_0)], & x\in I\setminus\{x_0\};\\
1, & x=x_0.
\end{dcases}
\end{equation}
\end{enumerate}
We call the function $T_n[f(x),x_0]$ of $x\in I\subseteq\mathbb{R}$ for $n\in\{-1\}\cup\mathbb{N}_0$ the $n$th normalized remainder of the formal series expansion~\eqref{Formal-series} of $f$ around $x_0$.
\par
When $x_0=0$, we simply write $T_n[f(x),0]$ as $T_n[f(x)]$.
\end{defn}

The terminology ``normalized remainder'' first explicitly appeared in the papers~\cite[Section~5]{log-exp-expan-Sym.tex}, \cite[Section~5]{OPENMATH-D-23-00112.tex}, \cite[Remark~7]{axioms-2962911.tex}, and~\cite{Niu-Qi-MDPI-Math.tex, llog-cosine-expan.tex, exp-norm-tail-ratio.tex, era-1998.tex, Filomat-24099.tex, log-secant-norm-tail.tex, JMI-4654.tex, Fractal-Fractional-Zhang-Yang-Qi-Du.tex}. Earlier literature also referred to this quantity as the ``normalized tail''; see~\cite[Remark~7]{axioms-2962911.tex} and~\cite{log-exp-expan-Sym.tex, Niu-Qi-MDPI-Math.tex, exp-norm-tail-ratio.tex, era-1998.tex, JMI-4654.tex, Fractal-Fractional-Zhang-Yang-Qi-Du.tex}. In addition, the term ``Qi's normalized remainder'' has been used in the papers~\cite{llog-cosine-expan.tex, log-secant-norm-tail.tex}, for example.
\par
The term ``normalized'' is used because $T_n[f(x_0),x_0]=1$ for all $n\ge-1$, distinguishing the quantity $T_n[f(x),x_0]$ from the classical remainders such as the Taylor, Lagrange, Cauchy, Schl\"omilch, and Rouch\'e remainders collected in~\cite[Section~6]{exp-tail-ratio.tex}. The notation $T_n[f(x)]$ was designed and applied in~\cite{Spivey-JIMA.tex, Spivey-ID-AP-JIMA.tex, ser-int-norm-tail.tex, very-short-exp-nt.tex, ratio-NR-squ-tan.tex, arcsine-liu-qi-mdpi.tex, JMI-5203.tex, JMI-5046.tex, exp-tail-ratio.tex, gam-recip-int-phi-sym.tex, Tan-Normalized-Remainder.tex}, for example.
\par
In the work, we will retrospect the origins of the concept of the normalized remainders of the Maclaurin power series expansions of functions, review main results on normalized remainders of some elementary functions since 2023, explore the historic backgrounds of the normalized remainders in combinatorics and number theory on the Bernoulli numbers and polynomials, the Stirling numbers and polynomials, and their generalizations by Broder, Carlitz, and Howard, present several basic properties of normalized remainders, and collect several conjectures and problems posed while investigating normalized remainders of several elementary functions.

\section{Origins in Mathematical Analysis}

We now start out to retrospect the origins of inventing the concept of the normalized remainders $T_n[f(x),x_0]$.

\subsection{First Origin Related to Tangent}\label{subsection2.1}
On 5 April 2023, Professor Chao-Ping Chen, a mathematician at Henan Polytechnic University, conjectured that the ratio $\frac{\ln F(x)}{\ln G(x)}$ is decreasing on $\bigl(0,\frac{\pi}{2}\bigr)$, where
\begin{equation}\label{f(x)-def-cases}
F(x)=\begin{dcases}
\frac{3(\tan x-x)}{x^3}, & 0<|x|<\frac{\pi}{2};\\
1, & x=0
\end{dcases}
\end{equation}
and
\begin{equation}\label{G(x)Tangent}
G(x)=\begin{dcases}
\frac{\tan x}x, & 0<|x|<\frac{\pi}{2};\\
1, & x=0.
\end{dcases}
\end{equation}
This conjecture was solved and applied in~\cite{NMCM-245722953.tex} and in the proof of~\cite[Theorem~3.1]{Cusa-type}.
\par
From the series expansion
\begin{align*}
\tan x&=\sum_{k=1}^{\infty}\frac{2^{2k}\bigl(2^{2k}-1\bigr)}{(2k)!}|B_{2k}|x^{2k-1}, \quad |x|<\frac{\pi}{2}\\
&=x+\frac{x^3}{3}+\frac{2 x^5}{15}+\frac{17 x^7}{315}+\frac{62 x^9}{2835}+\frac{1382 x^{11}}{155925}+\frac{21844 x^{13}}{6081075}+\dotsm
\end{align*}
in~\cite[p.~42]{Gradshteyn-Ryzhik-Table-8th}, we observed that
\begin{equation*}
F(x)=T_2[\tan x] \quad \text{and}\quad G(x)=T_0[\tan x].
\end{equation*}
\par
In~\cite{symmetry-2571630-english.tex}, we expanded the function $\ln T_2[\tan x]$ into the Maclaurin power series in two different forms. This result can be regarded as a generalization of the series expansion
\begin{align*}
\ln T_0[\tan x]&=\sum_{k=1}^{\infty}\frac{2^{2k}\bigl(2^{2k-1}-1\bigr)}{k(2k)!}|B_{2k}|x^{2k}\\
&=\frac{x^2}{3}+\frac{7 x^4}{90}+\frac{62 x^6}{2835}+\frac{127 x^8}{18900}+\frac{146 x^{10}}{66825}+\dotsm, \quad |x|<\frac{\pi}{2}
\end{align*}
in~\cite[p.~55]{Gradshteyn-Ryzhik-Table-8th}, where $B_{2k}$ denotes the Bernoulli numbers which are generated in~\cite[p.~3]{Temme-96-book} by
\begin{equation}\label{Bernoulli-GF}
\frac{z}{\te^z-1}=\sum_{k=0}^\infty B_k\frac{z^k}{k!}=1-\frac{z}2+\sum_{k=1}^\infty B_{2k}\frac{z^{2k}}{(2k)!}, \quad |z|<2\pi.
\end{equation}

\subsection{Second Origin Related to Sine}\label{2.2sec}
In the paper~\cite{log-sin-x-expan.tex}, motivated by the functions $F(x)$ and $G(x)$ defined in~\eqref{f(x)-def-cases} and~\eqref{G(x)Tangent}, we considered the functions
\begin{equation*}
H(x)=
\begin{dcases}
\frac{6}{x^2}\biggl(1-\frac{\sin x}{x}\biggr), & 0<|x|<\infty;\\
1, & x=0
\end{dcases}
\end{equation*}
and
\begin{equation}\label{R1(x)-x-3x-log}
R_1(x)=
\begin{dcases}
\frac{\ln H(x)}{\ln\frac{\sin x}{x}}, & |x|\in(0,\pi);\\
\frac{3}{10}, & x=0;\\
0, & x=\pm\pi.
\end{dcases}
\end{equation}
\par
Comparing with the well-known series expansion
\begin{align*}
\sin x&=\sum_{k=0}^{\infty}(-1)^k\frac{x^{2k+1}}{(2k+1)!}, \quad x\in\mathbb{R}\\
&=x-\frac{x^3}{6}+\frac{x^5}{120}-\frac{x^7}{5040}+\frac{x^9}{362880}-\dotsm,
\end{align*}
we observed that
\begin{equation}\label{T-sine-H(x)}
H(x)=T_2[\sin x] \quad\text{and}\quad \frac{\sin x}{x}=T_0[\sin x].
\end{equation}
\par
In the paper~\cite[Theorem~1]{log-sin-x-expan.tex}, motivated by the Maclaurin power series expansion of $\ln T_2[\tan x]$ obtained in~\cite{symmetry-2571630-english.tex}, we expanded the logarithm $\ln T_2[\sin x]$ into a Maclaurin series. This result can be regarded as a generalization of the series expansion
\begin{equation}\label{ln-sin-series}
\ln\frac{\sin x}{x}=-\sum_{k=1}^{\infty}\frac{2^{2k-1}}{(2k)!k}|B_{2k}|x^{2k}
=-\frac{x^2}{6}-\frac{x^4}{180}-\frac{x^6}{2835}-\dotsm
\end{equation}
for $|x|<\pi$ in~\cite[p.~55]{Gradshteyn-Ryzhik-Table-8th}.
\par
In the paper~\cite[Theorem~2]{log-sin-x-expan.tex}, stimulated by the decreasing monotonicity obtained in~\cite{NMCM-245722953.tex}, we proved that the ratio $R_1(x)$ defined by~\eqref{R1(x)-x-3x-log} decreasingly maps from $[0,\pi]$ onto $\bigl[0,\frac{3}{10}\bigr]$.

\subsection{Third Origin Related to Cosine}\label{2.3sec}
In the paper~\cite{OPENMATH-D-23-00112.tex}, motivated by the above observations and results in~\cite{symmetry-2571630-english.tex, log-sin-x-expan.tex, NMCM-245722953.tex}, we discussed the functions
\begin{equation*}
I(x)=
\begin{dcases}
\frac{2(1-\cos x)}{x^2}, & 0<|x|<2\pi;\\
1, & x=0
\end{dcases}
\end{equation*}
and
\begin{equation}\label{R2(x)-x-3x-log}
R_2(x)=
\begin{dcases}
\frac{\ln I(x)}{\ln\cos x}, & 0<|x|<\frac{\pi}{2};\\
\frac{1}{6}, & x=0;\\
0, & x=\pm\frac{\pi}{2}.
\end{dcases}
\end{equation}
\par
Comparing with the series expansion
\begin{align*}
\cos x&=\sum_{k=0}^{\infty}(-1)^k\frac{x^{2k}}{(2k)!}\\
&=1-\frac{x^2}{2}+\frac{x^4}{24}-\frac{x^6}{720}+\frac{x^8}{40320}-\dotsm, \quad x\in\mathbb{R},
\end{align*}
we observed that 
\begin{equation*}
I(x)=T_1[\cos x].
\end{equation*}
\par
In~\cite[Theorem~1]{OPENMATH-D-23-00112.tex}, we expanded the logarithm $\ln T_1[\cos x]$ into the Maclaurin series. This expansion can be regarded as a generalization of the series expansion
\begin{equation}
\begin{aligned}\label{log-cosine-series-expansion}
\ln\cos x&=-\sum_{k=1}^{\infty}\frac{2^{2k-1}(2^{2k}-1)}{k(2k)!}|B_{2k}|x^{2k}\\
&=-\frac{x^2}{2}-\frac{x^4}{12}-\frac{x^6}{45}-\frac{17x^8}{2520}-\dotsm, \quad |x|<\frac{\pi}2,
\end{aligned}
\end{equation}
which can be found in~\cite[p.~42]{Gradshteyn-Ryzhik-Table-8th}.
\par
In~\cite[Theorem~2]{OPENMATH-D-23-00112.tex}, we proved that the ratio $R_2(x)$ defined by~\eqref{R2(x)-x-3x-log} decreasingly maps $\bigl[0,\frac{\pi}2\bigr]$ onto $\bigl[0,\frac{1}{6}\bigr]$.
\par
In~\cite[Remark~2]{OPENMATH-D-23-00112.tex}, Theorem~1 in~\cite{OPENMATH-D-23-00112.tex} was simplified as
\begin{equation}\label{F(x)-ser-Bernoulli-No}
\ln T_1[\cos x]=-\sum_{k=1}^{\infty}\frac{|B_{2k}|}{k}\frac{x^{2k}}{(2k)!}, \quad |x|<2\pi.
\end{equation}
In~\cite[Remark~3]{OPENMATH-D-23-00112.tex}, we constructed a positive and even function
\begin{equation}\label{Jn(x)-eq}
J_n(x)=
\begin{dcases}
-\frac{n+1}{|B_{2n+2}|}\frac{(2n+2)!}{x^{2n+2}} \Biggl[\ln T_1[\cos x]+\sum_{k=1}^{n}\frac{|B_{2k}|}{k}\frac{x^{2k}}{(2k)!}\Biggr], & 0<|x|<2\pi\\
1, & x=0
\end{dcases}
\end{equation}
for $n\in\mathbb{N}_0$. We observed that
\begin{equation*}
J_n(x)=T_{2n+1}[\ln T_1[\cos x]], \quad n\in\mathbb{N}_0.
\end{equation*}
\par
Further, in the conclusion section~\cite[Section~5]{OPENMATH-D-23-00112.tex}, we first explicitly invented the concept and used the terminology of the normalized remainder as follows.
\begin{quote}
Let $f(x)$ be an even, positive, and analytic function on $(-r,r)$ such that $f(0)=1$ and $f^{(2m)}(0)\ne0$ for $m\ge1$. Then
\begin{equation*}
f(x)=\sum_{k=0}^{\infty}f^{(2k)}(0)\frac{x^{2k}}{(2k)!}, \quad x\in(-r,r).
\end{equation*}
What are the Maclaurin power series expansions of the logarithmic expressions $\ln f(x)$ and $\ln\frac{2[f(x)-1]}{f''(0)x^2}$? What is the monotonicity of the even function $\frac{\ln\frac{2[f(x)-1]}{f''(0)x^2}}{\ln f(x)}$ on the interval $(0,r)$? Generally, what about the properties for the logarithms of the normalized remainders
\begin{equation*}
F_n(x)=
\begin{dcases}
0, & x=0\\
\ln\Biggl\{\frac{(2n)!}{f^{(2n)}(0)}\frac1{x^{2n}}\Biggl[f(x)-\sum_{k=0}^{n-1}f^{(2k)}(0)\frac{x^{2k}}{(2k)!}\Biggr]\Biggr\}, & x\ne0
\end{dcases}
\end{equation*}
for $n\ge0$ and their ratios $R_{m,n}(x)=\frac{F_n(x)}{F_m(x)}$ for $n>m\ge0$?
\end{quote}

\section{Recent Results of Normalized Remainders}
We now in a position to review main results of normalized remainders of several elementary functions.

\subsection{Normalized Remainders of Sine and Cosine Functions}

As a continuation of Sections~\ref{2.2sec} and~\ref{2.3sec}, we considered the general normalized remainders of the cosine and sine functions in several works.

\subsubsection{Normalized Remainder of Cosine}
In the paper~\cite{era-1998.tex}, we considered the normalized remainder
\begin{equation*}
T_{2n-1}[\cos x]=
\begin{dcases}
(-1)^n\frac{(2n)!}{x^{2n}}\Biggl[\cos x-\sum_{k=0}^{n-1}(-1)^k\frac{x^{2k}}{(2k)!}\Biggr], & x\ne0\\
1, & x=0
\end{dcases}
\end{equation*}
for $n\in\mathbb{N}_0$, where an empty sum is understood to be $0$.
\par
It is easy to see that
\begin{equation*}
T_{-1}[\cos x]=\cos x \quad\text{and}\quad T_1[\cos x]=I(x).
\end{equation*}
\par
In~\cite[Theorem~1]{era-1998.tex}, we expanded the logarithm $\ln T_{2n-1}[\cos x]$ for $n\in\mathbb{N}_0$ into a Maclaurin series. This gives a generalization of~\cite[Theorem~1]{OPENMATH-D-23-00112.tex} and the series expansions~\eqref{log-cosine-series-expansion} and~\eqref{F(x)-ser-Bernoulli-No}.
\par
In~\cite[Theorem~2]{era-1998.tex}, we proved the following conclusions:
\begin{enumerate}[1)]
\item
For $n=0$ and $n\ge2$, the normalized remainder $T_{2n-1}[\cos x]$ is decreasing and logarithmically concave on $\bigl(0,\frac{\pi}{2}\bigr)$. 
\item
The normalized remainder $T_1[\cos x]$ is decreasing on the intervals $(0,2\pi)$ and $(x_k,2(k+1)\pi)$ for $k\in\mathbb{N}$, and increasing on $(2k\pi,x_k)$ for $k\in\mathbb{N}$, where $x_k\in(2k\pi,2(k+1)\pi)$ for $k\in\mathbb{N}$ is the zero of the equation $\tan\frac{x}{2}=\frac{x}{2}$ on $(0,\infty)$.
\end{enumerate}
\par
In~\cite[Theorem~3]{era-1998.tex}, we verified that the ratio
\begin{equation*}
R_3(x)=
\begin{dcases}
\frac{\ln T_3[\cos x]}{\ln\cos x}, &0<|x|<\frac{\pi}{2}\\
\frac{1}{15}, & x=0\\
0, & x=\pm\frac{\pi}{2}
\end{dcases}
\end{equation*}
is decreasing on $\bigl[0,\frac{\pi}2\bigr]$. This monotonicity result generalizes~\cite[Theorem~2]{OPENMATH-D-23-00112.tex}, because the ratio $R_3(x)$ is a generalization of $R_2(x)$ defined in~\eqref{R2(x)-x-3x-log}.

\subsubsection{Normalized Remainder of Sine}
In~\cite{Fractal-Fractional-Zhang-Yang-Qi-Du.tex}, we discussed the normalized remainder $T_{2n-1}[\cos x]$ and
\begin{equation*}
T_{2n}[\sin x]=
\begin{dcases}
(-1)^n\frac{(2n+1)!}{x^{2n+1}}\Biggl[\sin x-\sum_{k=0}^{n-1} (-1)^k\frac{x^{2k+1}}{(2k+1)!}\Biggr], & x\ne0\\
1, & x=0
\end{dcases}
\end{equation*}
for $n\in\mathbb{N}_0$, where an empty sum is understood to be $0$.
\par
It is obvious that those equalities in~\eqref{T-sine-H(x)} are valid.
\par
Among other things, the main results in~\cite{Fractal-Fractional-Zhang-Yang-Qi-Du.tex} can be recited as follows:
\begin{enumerate}[1)]
\item
The normalized remainders $T_{2n}[\sin x]$ and $T_{2n-1}[\cos x]$ for $n\in\mathbb{N}$ and $x\in\mathbb{R}$ have the integral representations
\begin{align*}
T_{2n}[\sin x]&=
\begin{dcases}
2n(2n+1)\int_0^1 (1-u)^{2n-1}\frac{\sin(xu)}{x}\td u, & x\ne0\\
1, & x=0
\end{dcases}\\
&=(2n+1)\int_0^1 (1-u)^{2n}\cos(xu)\td u
\end{align*}
and
\begin{align*}
T_{2n-1}[\cos x]&=
\begin{dcases}
2n(2n-1)\int_0^1 (1-u)^{2n-2} \frac{\sin(xu)}{x}\td u, & x\ne0\\
1, & x=0
\end{dcases}\\
&=2n\int_0^1 (1-u)^{2n-1} \cos(xu)\td u.
\end{align*}
\item
For $n\in\mathbb{N}$, the normalized remainder $T_{2n}[\sin x]$ is positive and decreasing in $x\in(0,\infty)$, and concave in $x\in(0,\pi)$.
\item
The normalized remainder $T_1[\cos x]$ is nonnegative on $(0,\infty)$, decreasing on $[0,2\pi]$, and concave on $(0,x_0)$, where $x_0\in\bigl(\frac{\pi}{2},\pi\bigr)$ is the first positive zero of the equation
\begin{equation*}
\bigl(x^2-2\bigr) \sin x+2x\cos x=0.
\end{equation*}
For $n\ge2$, the normalized remainder $T_{2n-1}[\cos x]$ is positive and decreasing on $(0,\infty)$, and concave on $(0,\pi)$.
\item
For $n\in\mathbb{N}$, the normalized remainders $T_{2n}[\sin x]$ and $T_{2n-1}[\cos x]$ have the relations
\begin{equation}\label{1F2-Sin-Special=Value}
T_{2n}[\sin x]={}_1F_2\biggl(1; n+1, n+\frac{3}{2}; -\frac{x^2}{4}\biggr)
\end{equation}
and
\begin{equation}\label{1F2-Cos-Special=Value}
T_{2n-1}[\cos x]={}_1F_2\biggl(1;n+\frac{1}{2}, n+1; -\frac{x^2}{4}\biggr),
\end{equation}
where the generalized hypergeometric function ${}_1F_2(a;b,c;z)$ is defined in~\cite[p.~1020]{Gradshteyn-Ryzhik-Table-8th} by
\begin{equation*}
{}_1F_2(a;b,c;z)=\sum_{n=0}^\infty\frac{(a)_n} {(b)_n(c)_n}\frac{z^n}{n!},\quad z\in\mathbb{C}
\end{equation*}
for $a\in\mathbb{C}$, $b,c\in\mathbb{C}\setminus\{0,-1,-2,\dotsc\}$, and
\begin{equation*}
(z)_n=\prod_{\ell=0}^{n-1}(z+\ell)
=
\begin{cases}
z(z+1)\dotsm(z+n-1), & n\in\mathbb{N};\\
1, & n=0.
\end{cases}
\end{equation*}
See also~\cite[Remark~7]{axioms-2962911.tex}.
\end{enumerate}

\subsubsection{Ratios of Normalized Remainders of Sine and Cosine}
In the paper~\cite{Niu-Qi-MDPI-Math.tex}, we considered the ratios
\begin{equation*}
\frac{T_{2n-1}[\cos x]}{T_{2n}[\sin x]}, \quad \frac{T_{2n-1}[\cos x]}{T_{2n+1}[\cos x]}, \quad \frac{T_{2n}[\sin x]}{T_{2n+2}[\sin x]}
\end{equation*}
in $x\in(0,\infty)$ for $n\in\mathbb{N}$. We obtained the following conclusions:
\begin{enumerate}[1)]
\item
The limits
\begin{align*}
\lim_{x\to0}\frac{T_{2n-1}[\cos x]}{T_{2n}[\sin x]}&=1, & \lim_{x\to0}\frac{T_{2n-1}[\cos x]}{T_{2n+1}[\cos x]}&=1, \\
\lim_{x\to0}\frac{T_{2n}[\sin x]}{T_{2n+2}[\sin x]}&=1, & \lim_{x\to\infty}\frac{T_{2n}[\sin x]}{T_{2n+2}[\sin x]}&=\frac{n(2n+1)}{(n+1)(2n+3)}
\end{align*}
are valid for $n\in\mathbb{N}$, while the limits
\begin{equation*}
\lim_{x\to\infty}\frac{T_{2n-1}[\cos x]}{T_{2n}[\sin x]}=\frac{2n-1}{2n+1}
\end{equation*}
and
\begin{equation*}
\lim_{x\to\infty}\frac{T_{2n-1}[\cos x]}{T_{2n+1}[\cos x]}=\frac{n(2n-1)}{(n+1)(2n+1)}
\end{equation*}
are valid for $n\ge2$.
\item
For $n\in\mathbb{N}$, the limit
\begin{equation*}
\lim_{x\to\infty}\frac{\int_0^1 (1-u)^{2n}\cos(xu)\td u}{\int_0^1 (1-u)^{2n+2}\cos(xu)\td u}=\frac{n}{n+1}
\end{equation*}
is valid. For $n\ge2$, the limits
\begin{equation*}
\lim_{x\to\infty}\frac{\int_0^1 (1-u)^{2n-1} \cos(xu)\td u} {\int_0^1 (1-u)^{2n}\cos(xu)\td u}=\frac{2n-1}{2n}
\end{equation*}
and
\begin{equation*}
\lim_{x\to\infty}\frac{\int_0^1 (1-u)^{2n-1} \cos(xu)\td u} {\int_0^1 (1-u)^{2n+1} \cos(xu)\td u}=\frac{2n-1}{2n+1}
\end{equation*}
are valid.
\item
The ratios $\frac{T_1[\cos x]}{T_2[\sin x]}$ and $\frac{T_1[\cos x]}{T_3[\cos x]}$ have infinitely many minimums $0$ at $x=2k\pi$ for $k\in\mathbb{N}$, the ratios
\begin{align*}
\frac{T_1[\cos x]}{T_2[\sin x]}, \quad
\frac{T_1[\cos x]}{T_3[\cos x]}, \quad
\frac{T_2[\sin x]}{T_4[\sin x]}, \quad
\frac{T_3[\cos x]}{T_4[\sin x]}, \quad
\frac{T_3[\cos x]}{T_5[\cos x]}
\end{align*}
are not monotonic in $x\in(0,\infty)$, and the ratio $\frac{T_4[\sin x]}{T_6[\sin x]}$ is decreasing in $x\in(0,\infty)$.
\item
For $n\ge3$, the ratios
\begin{equation*}
\frac{T_{2n-1}[\cos x]}{T_{2n}[\sin x]}, \quad \frac{T_{2n-1}[\cos x]}{T_{2n+1}[\cos x]}, \quad \frac{T_{2n}[\sin x]}{T_{2n+2}[\sin x]}
\end{equation*}
are decreasing in $x\in(0,\infty)$.
\item
For $x\in(0,\infty)$, the following double inequalities hold:
\begin{align*}
\frac{2n-1}{2n+1}&<\frac{T_{2n-1}[\cos x]}{T_{2n}[\sin x]}<1, \quad n\ge3,\\
\frac{n(2n-1)}{(n+1)(2n+1)}&<\frac{T_{2n-1}[\cos x]}{T_{2n+1}[\cos x]}<1,\quad n\ge3,
\intertext{and}
\frac{n(2n+1)}{(n+1)(2n+3)}&<\frac{T_{2n}[\sin x]}{T_{2n+2}[\sin x]}<1, \quad n\ge2.
\end{align*}
\item
Using the formulas~\eqref{1F2-Sin-Special=Value} and~\eqref{1F2-Cos-Special=Value}, as done in~\cite[Remark~7]{axioms-2962911.tex}, we reformulated in~\cite[Section~5]{Niu-Qi-MDPI-Math.tex} all the above results related to the normalized remainders $T_{2n}[\sin x]$ and $T_{2n-1}[\cos x]$ in terms of the generalized hypergeometric function ${}_1F_2(a;b,c;z)$.
\end{enumerate}

\subsubsection{Normalized Remainder of Logarithm of Normalized Remainder for Cosine}
In the paper~\cite{llog-cosine-expan.tex}, the authors investigated the normalized remainder 
\begin{equation*}
J_n(x)=T_{2n+1}[\ln T_1[\cos x]]
\end{equation*} 
defined by~\eqref{Jn(x)-eq} for $n\in\mathbb{N}_0$. They discovered the following results:
\begin{enumerate}[1)]
\item
For $n\in\mathbb{N}_0$, the normalized remainder $T_{2n+1}[\ln T_1[\cos x]]$ is
\begin{enumerate}
\item
even in $x\in(-2\pi,2\pi)$,
\item
increasing in $x\in(0,2\pi)$,
\item
logarithmically convex in $x\in(-2\pi,2\pi)$.
\end{enumerate}
\item
For $n\in\mathbb{N}_0$, the logarithm $\ln T_{2n+1}[\ln T_1[\cos x]]$ was expanded into a Maclaurin series.
\item
For $n\in\mathbb{N}_0$, the ratio $\frac{T_{2n+3}[\ln T_1[\cos x]]}{T_{2n+1}[\ln T_1[\cos x]]}$ is increasing in $x\in(0,2\pi)$.
\end{enumerate}

\subsubsection{Normalized Remainder of Logarithm of Secant}
In the paper~\cite{log-secant-norm-tail.tex}, the authors considered the normalized remainder
\begin{multline*}
T_{2n+1}[\ln\sec x]=\\
\begin{dcases}
\frac{(n+1)(2n+2)!\Bigl[\ln\sec x-\sum_{k=1}^{n}\frac{2^{2k-1}(2^{2k}-1)}{k(2k)!}|B_{2k}|x^{2k}\Bigr]} {2^{2n+1}(2^{2n+2}-1)|B_{2n+2}|x^{2n+2}}, & x\ne0\\
1, & x=0
\end{dcases}
\end{multline*}
for $n\in\mathbb{N}_0$ and $x\in\bigl(-\frac{\pi}{2},\frac{\pi}{2}\bigr)$, where an empty sum is understood to be $0$.
\par
Since $\ln\sec x=-\ln\cos x$, it easily follows that
\begin{equation}\label{ln-sec-cos-tail}
T_{2n+1}[\ln\sec x]=T_{2n+1}[\ln\cos x]
\end{equation}
for $n\in\mathbb{N}_0$ and $x\in\bigl(-\frac{\pi}{2},\frac{\pi}{2}\bigr)$.
\par
In the paper~\cite{log-secant-norm-tail.tex}, they gained the following two conclusions:
\begin{enumerate}[1)]
\item
The normalized remainder $T_{2n+1}[\ln\sec x]$ for given $n\in\mathbb{N}_0$ is a logarithmically convex function in $x\in\bigl(-\frac{\pi}{2},\frac{\pi}{2}\bigr)$.
\item
The ratio $\frac{T_{2n+3}[\ln\sec x]}{T_{2n+1}[\ln\sec x]}$ for $n\in\mathbb{N}_0$ is an increasing function of $x\in\bigl(0,\frac{\pi}{2}\bigr)$ and a decreasing function of  $x\in\bigl(-\frac{\pi}{2},0\bigr)$.
\end{enumerate}

\subsection{Normalized Remainders of Tangent Function and its Square}
Recall from the monographs~\cite[Chapter~XIII]{mpf-1993}, \cite{Schilling-Song-Vondracek-2nd}, and~\cite[Chapter~IV]{Widder-Laplace-Transform-41} that,
\begin{enumerate}[1)]
\item
an infinitely differentiable real function $f(x)$ defined on an interval $I$ is said to be absolutely monotonic in $x\in I$ if and only if all of its derivatives satisfy $f^{(n)}(x)\ge0$ for $n\in\mathbb{N}_0$ and $x\in I$;
\item
an infinitely differentiable real function $f(x)$ defined on $I$ is said to be completely monotonic in $x\in I$ if and only if all of its derivatives satisfy $(-1)^nf^{(n)}(x)\ge0$ for $n\in\mathbb{N}_0$ and $x\in I$.
\end{enumerate}
A function $f(x)$ is completely monotonic on $(a,b)$ if and only if it is absolutely monotonic on $(-b,-a)$; see~\cite[p.~145, Definition~2c]{Widder-Laplace-Transform-41}.
\par
Recall from~\cite[Definition~1]{absolute-mon-simp.tex} and~\cite[Definition~1]{compmon2} that,
\begin{enumerate}[1)]
\item
a positive function $f(x)$ is said to be logarithmically absolutely monotonic on an interval $I$ if it has derivatives of all orders and $[\ln f(x)]^{(n)}\ge0$ for $x\in I$ and $n\in\mathbb{N}$.
\item
a positive function $f(x)$ is said to be logarithmically completely monotonic on an interval $I$ if it has derivatives of all orders and the inequality $(-1)^n[\ln f(x)]^{(n)}\ge0$ holds for $x\in I$ and $n\in\mathbb{N}$.
\end{enumerate}
\par
In~\cite[Theorem~1]{absolute-mon-simp.tex}, the authors proved that a logarithmically absolutely monotonic function on an interval $I$ is also absolutely monotonic on $I$, but not conversely. In the paper~\cite{CBerg}, \cite[Theorem~4]{absolute-mon-simp.tex}, and~\cite[Theorem~1]{compmon2}, the authors proved that a logarithmically completely monotonic function on an interval $I$ is also completely monotonic on $I$, but not conversely.

\subsubsection{Normalized Remainders of Tangent}
As a continuation of the study in Section~\ref{subsection2.1}, we considered in the paper~\cite{Tan-Normalized-Remainder.tex} the general normalized remainder
\begin{equation*}
T_{2n-2}[\tan x]=
\begin{dcases}
\frac{(2n)!\bigl[\tan x-\sum_{k=1}^{n-1}\frac{2^{2k}(2^{2k}-1)}{(2k)!}|B_{2k}|x^{2k-1}\bigr]} {2^{2n}\bigl(2^{2n}-1\bigr)|B_{2n}|x^{2n-1}}, & 0<|x|<\frac{\pi}{2}\\
1, & x=0
\end{dcases}
\end{equation*}
for $n\in\mathbb{N}$. Our main results are the following ones:
\begin{enumerate}[1)]
\item
For $n\in\mathbb{N}$, the normalized remainder $T_{2n-2}[\tan x]$ is even on $\bigl(-\frac{\pi}{2},\frac{\pi}{2}\bigr)$, absolutely monotonic in $x\in\bigl(0,\frac{\pi}{2}\bigr)$, and completely monotonic in $x\in\bigl(-\frac{\pi}{2},0\bigr)$.
\item
For $n\in\mathbb{N}$, the normalized remainder $T_{2n-2}[\tan x]$ is logarithmically absolutely monotonic in $x\in\bigl(0,\frac{\pi}{2}\bigr)$ and logarithmically completely monotonic in $x\in\bigl(-\frac{\pi}{2},0\bigr)$.
\item
For $n\in\mathbb{N}$, the ratio $\frac{(T_{2n-2}[\tan x])'}{x T_{2n-2}[\tan x]}$ is absolutely monotonic in $x\in\bigl(0,\frac{\pi}{2}\bigr)$ and completely monotonic in $x\in\bigl(-\frac{\pi}{2},0\bigr)$.
\item
For $n\in\mathbb{N}$, the function $1-\frac{T_{2n-2}[\tan x]}{T_{2n}[\tan x]}$ is absolutely monotonic in $x\in\bigl(0,\frac{\pi}{2}\bigr)$ and completely monotonic in $x\in\bigl(-\frac{\pi}{2},0\bigr)$.
Consequently, the ratio $\frac{T_{2n-2}[\tan x]}{T_{2n}[\tan x]}$ for $n\in\mathbb{N}$ is decreasing in $x\in\bigl(0,\frac{\pi}{2}\bigr)$, increasing in $x\in\bigl(-\frac{\pi}{2},0\bigr)$, and concave in $x\in\bigl(-\frac{\pi}{2},\frac{\pi}{2}\bigr)$.
\item
For $n\in\mathbb{N}$, the logarithm $\ln T_{2n-2}[\tan x]$ was expanded into a Maclaurin power series.
\end{enumerate}

\subsubsection{Normalized Remainders of Squared Tangent}
In~\cite{JMI-4654.tex}, we built the normalized remainder
\begin{multline*}
T_{2n-1}[\tan^2x]=\\
\begin{dcases}
\frac{\displaystyle(2n+2)!\Biggl[\tan^2x-\sum_{\ell=1}^{n-1} \frac{2^{2\ell+2}\bigl(2^{2\ell+2}-1\bigr)(2\ell+1)}{(2\ell+2)!} |B_{2\ell+2}|x^{2\ell}\Biggr]} {2^{2n+2}(2^{2n+2}-1)(2n+1)|B_{2n+2}|x^{2n}}, & x\ne0\\
1, & x=0
\end{dcases}
\end{multline*}
$n\in\mathbb{N}$ and $|x|<\frac{\pi}{2}$.
\par
Due to the identity $\sec^2x=1+\tan^2x$, we derived in~\cite{JMI-4654.tex} that
\begin{equation}\label{tan-sec-square}
T_{2n-1}[\tan^2x]=T_{2n-1}[\sec^2x], \quad n\in\mathbb{N}.
\end{equation}
\par
Our main results in~\cite{JMI-4654.tex} are as follows:
\begin{enumerate}[1)]
\item
For $n\in\mathbb{N}$, the normalized remainder $T_{2n-1}[\tan^2x]$ is even, positive, and increasing on $\bigl(0,\frac{\pi}{2}\bigr)$.
\item
For $n\in\mathbb{N}$, the normalized remainder $T_{2n-1}[\tan^2x]$ is logarithmically convex on $\bigl(-\frac{\pi}{2},\frac{\pi}{2}\bigr)$.
\item
For $n\in\mathbb{N}$, the logarithm $\ln T_{2n-1}[\tan^2x]$ was expanded into a Maclaurin series.
\end{enumerate}

\subsubsection{Ratio of Normalized Remainders of Squared Tangent}
In~\cite{ratio-NR-squ-tan.tex}, the ratio $\frac{T_{2n+1}[\tan^2x]}{T_{2n-1}[\tan^2x]}$ for $n\in\mathbb{N}$ was proved to be increasing in $x\in\bigl(-\frac{\pi}{2},0\bigr)$ and decreasing in $x\in\bigl(0,\frac{\pi}{2}\bigr)$.

\subsection{Normalized Remainders of Exponential Function}
\subsubsection{Generating Function for Bernoulli Numbers}\label{sec3.31}
In the paper~\cite{Bernouli-No-Tail.tex}, we considered the function $\frac{x}{\te^x-1}$ in~\eqref{Bernoulli-GF} and studied its normalized remainder
\begin{equation*}
T_{2n+1}\biggl[\frac{x}{\te^x-1}\biggr]=
\begin{dcases}
\frac{(2n+2)!}{B_{2n+2}} \frac{1}{x^{2n+2}} \Biggl[\frac{x}{\te^x-1}-1+\frac{x}2-\sum_{j=1}^n B_{2j}\frac{x^{2j}}{(2j)!}\Biggr], & x\ne0\\
1, & x=0
\end{dcases}
\end{equation*}
for $n\in\mathbb{N}_0$ and $x\in\mathbb{R}$, where an empty sum is assumed to be $0$.
\par
Our main results in~\cite{Bernouli-No-Tail.tex} are recited as follows:
\begin{enumerate}[1)]
\item
For $n\in\mathbb{N}$, the normalized remainder $T_{2n+1}\bigl[\frac{x}{\te^x-1}\bigr]$ is even on $(-\infty,\infty)$.
\item
For $n\in\mathbb{N}$, the normalized remainder $T_{2n+1}\bigl[\frac{x}{\te^x-1}\bigr]$ is positive in $x\in(0,\infty)$.
\item
For $n\in\mathbb{N}$, the normalized remainder $T_{2n+1}\bigl[\frac{x}{\te^x-1}\bigr]$ is decreasing in $x\in(0,\infty)$.
\item
For $n\in\mathbb{N}$, the ratio $\frac{T_{2n+1}[\frac{x}{\te^x-1}]}{T_{2n+3}[\frac{x}{\te^x-1}]}$ is increasing in $x\in(0,\infty)$.
\item
For $n\in\mathbb{N}$, the ratio $\frac{B_{2n-1}(t)}{B_{2n+1}(t)}$ is increasing in $t\in\bigl(0,\frac{1}{2}\bigr)$ and decreasing in $t\in\bigl(\frac{1}{2},1\bigr)$, where the classical Bernoulli polynomials $B_k(t)$ are generated in~\cite[p.~3]{Temme-96-book} by
\begin{equation}\label{Bernoulli-Polyn-GF}
\frac{x\te^{x t}}{\te^x-1}=\sum_{k=0}^{\infty}B_k(t)\frac{x^k}{k!}, \quad |x|<2\pi.
\end{equation}
\end{enumerate}
\par
We noticed that the generating function $\frac{x}{\te^x-1}$ of the Bernoulli numbers $B_k$ is just the reciprocal $\frac1{T_0[\te^x]}$ and that the generating function $\frac{x\te^{x t}}{\te^x-1}$ of the Bernoulli polynomials $B_k(t)$ is just the reciprocal $\frac1{T_0[(\te^x-1)\te^{-t x}]}$ or $\frac{\te^{x t}}{T_0[\te^x]}$. As a result, we derive
\begin{equation*}
T_{2n+1}\biggl[\frac{x}{\te^x-1}\biggr]=T_{2n+1}\biggl[\frac1{T_0[\te^x]}\biggr]
\end{equation*}
for $n\in\mathbb{N}_0$ and $x\in\mathbb{R}$.
\par
By~\cite[Theorem~2.1]{Guo-Qi-Filomat-2011-May-12.tex}, we deduce that the logarithm $\ln T_0[\te^x]$ is convex on $(-\infty,\infty)$, $3$-convex (that is, $\frac{\operatorname{d}^3\ln T_0[\te^x]}{\td x^3}\ge0$) on $(-\infty,0)$, and $3$-concave (that is, $\frac{\operatorname{d}^3\ln T_0[\te^x]}{\td x^3}\le0$) on $(0,\infty)$.
\par
Motivated by~\cite[Proposition~1]{Bernouli-No-Tail.tex}, Yang--Qi~\cite{arxiv.2405.05280} and Pinelis~\cite{Pinelis-Lie-MIA-24} further investigated the monotonicity and other properties of the ratios
\begin{equation*}
\frac{B_{2\ell-1}(s)}{B_{2\ell+1}(s)}, \quad \frac{B_{2\ell}(s)}{B_{2\ell+1}(s)},\quad \frac{B_{2m}(s)}{B_{2\ell}(s)},\quad \frac{B_{2\ell}(s)}{B_{2\ell-1}(s)},\quad \ell,m\in\mathbb{N}.
\end{equation*}
\par
Since 2023, when studying properties of normalized remainders of the Maclaurin power series expansions, we discovered some monotonicity results and inequalities on sequences and functions involving the ratio $\bigl|\frac{B_{2k+2}}{B_{2k}}\bigr|$ and the Riemann zeta function $\zeta(t)$. For more details, please refer to the papers~\cite{MIA-9509.tex} and~\cite[Section~1.12]{arcsine-liu-qi-mdpi.tex}.

\subsubsection{Normalized Remainder of Exponential Function}
In the paper~\cite{log-exp-expan-Sym.tex}, we obtained the following results:
\begin{enumerate}[1)]
\item
The logarithm $\ln T_0[\te^x]$ was expanded into the series
\begin{align*}
\ln T_0[\te^x]&=\frac{x}{2}+\sum_{k=1}^\infty\frac{B_{2k}}{2k}\frac{x^{2k}}{(2k)!}\\
&=\frac{x}{2}+\frac{x^2}{24}-\frac{x^4}{2880}+\frac{x^6}{181440} -\frac{x^8}{9676800}+\dotsm, \quad |x|<2\pi.
\end{align*}
\item
The logarithm $\ln T_1[\te^x]$ was expanded into the series
\begin{align*}
\ln T_1[\te^x]&=\frac{x}{3}+2\sum_{k=1}^{\infty}\frac{A_{k+1}}{k+1}\frac{x^{k+1}}{(k+1)!}\\
&=\frac{x}{3}+\frac{x^2}{36}+\frac{x^3}{810}-\frac{x^4}{12960}-\frac{x^5}{68040}-\frac{x^6}{12247200}+\frac{x^7}{6123600}+\dotsm,
\end{align*}
where $A_k$ are generated (see~\cite[p.~979, Eq.~(2.9)]{Howard-MC-80-977}) by
\begin{equation}
\begin{aligned}\label{Howard-numbers-def}
\frac{1}{T_1[\te^x]}&=\sum_{k=0}^{\infty}A_k\frac{x^k}{k!}\\
&=1-\frac{x}{3}+\frac{x^2}{36}+\frac{x^3}{540}-\frac{x^4}{6480}-\frac{x^5}{27216}-\frac{x^6}{4082400}+\dotsm
\end{aligned}
\end{equation}
for $|x|<|x_0|$, where $x_0\ne0$ is the zero, closest to the origin $x=0$, of the equation $\te^x-1-x=0$ on the complex plane $\mathbb{C}$.
\par
In~\cite[Theorem~2.1]{Bell-value-elem-funct.tex}, a closed-form expression for $A_j$ was given by
\begin{multline}\label{A(k)Bell-Eq}
A_k=\frac{1}{(2k)!!}\sum_{\ell=0}^k \frac{\binom{k}{\ell}}{\binom{k+\ell}{\ell}} \sum_{n=0}^{k-\ell}n!\langle-2\ell\rangle_{k-\ell-n} \binom{k-\ell}{n}\binom{k+\ell}{n} \\
\times\sum_{m=0}^{\ell}(-1)^{m} \binom{\ell}{m} \sum_{p=0}^{k+\ell} 2^p\langle \ell-m\rangle_{k+\ell-p}\binom{k+\ell}{p} \frac{S(m+p,m)}{\binom{m+p}{m}},
\end{multline}
where the falling factorial $\langle\rho\rangle_k$ is defined by
\begin{equation*}
\langle\rho\rangle_k=
\prod_{j=0}^{k-1}(\rho-j)=
\begin{cases}
\rho(\rho-1)\dotsm(\rho-k+1), & k\in\mathbb{N}\\
1,& k=0
\end{cases}
\end{equation*}
for $\rho\in\mathbb{C}$ and the Stirling numbers of the second kind $S(j,k)$ for $j\ge k\in\mathbb{N}_0$ can be analytically generated~\cite[p.~51]{Comtet-Combinatorics-74} by
\begin{equation}\label{2Stirl-funct-rew}
(T_0[\te^x])^\ell=\sum_{k=0}^\infty \frac{S(k+\ell,\ell)}{\binom{k+\ell}{\ell}} \frac{x^{k}}{k!}, \quad \ell\in\mathbb{N}_0.
\end{equation}
\item
The logarithm $\ln T_{n-1}[\te^x]$ of the normalized remainder
\begin{equation*}
T_{n-1}[\te^x]=
\begin{dcases}
\frac{n!}{x^{n}}\Biggl(\te^x-\sum_{k=0}^{n-1}\frac{x^k}{k!}\Biggr), & x\ne0\\
1, & x=0
\end{dcases}
\end{equation*}
for $n\in\mathbb{N}$ and $x\in\mathbb{R}$ was expanded into a Maclaurin series.
\item
The Bernoulli numbers $B_{2k}$ and the Howard numbers $A_k$ were expressed by $B_{2k}=-2kD_{2k}(1)$ for $k\in\mathbb{N}$ and by
\begin{equation}\label{Ak=determ}
A_k=(-1)^{k-1}\frac{k}{2}D_{k}(2), \quad k\ge2,
\end{equation}
where the determinant $D_{\ell}(n)$ is defined by
\begin{align*}
D_{\ell+1}(n)
&=\begin{vmatrix}
\frac{1}{\binom{n+1}{1}} & \frac{\binom{0}{0}}{\binom{n}{0}} & 0 & 0 & \dotsm & 0 & 0\\
\frac{1}{\binom{n+2}{2}} & \frac{\binom{1}{0}}{\binom{n+1}{1}} & \frac{\binom{1}{1}}{\binom{n}{0}} & 0 & \dotsm & 0 & 0\\
\frac{1}{\binom{n+3}{3}} & \frac{\binom{2}{0}}{\binom{n+2}{2}} & \frac{\binom{2}{1}}{\binom{n+1}{1}} & \frac{\binom{2}{2}}{\binom{n}{0}} & \dotsm & 0  & 0\\
\frac{1}{\binom{n+4}{4}} & \frac{\binom{3}{0}}{\binom{n+3}{3}} & \frac{\binom{3}{1}}{\binom{n+2}{2}} & \frac{\binom{3}{2}}{\binom{n+1}{1}} & \dotsm & 0  & 0\\
\frac{1}{\binom{n+5}{5}} & \frac{\binom{4}{0}}{\binom{n+4}{4}} & \frac{\binom{4}{1}}{\binom{n+3}{3}} & \frac{\binom{4}{2}}{\binom{n+2}{2}} & \dotsm & 0  & 0\\
\vdots & \vdots & \vdots & \vdots & \ddots & \vdots & \vdots\\
\frac{1}{\binom{n+\ell-1}{\ell-1}} & \frac{\binom{\ell-2}{0}}{\binom{n+\ell-2}{\ell-2}} & \frac{\binom{\ell-2}{1}}{\binom{n+\ell-3}{\ell-3}} & \frac{\binom{\ell-2}{2}}{\binom{n+\ell-4}{\ell-4}} & \dotsm & \frac{\binom{\ell-2}{\ell-2}}{\binom{n}{0}} & 0\\
\frac{1}{\binom{n+\ell}{\ell}} & \frac{\binom{\ell-1}{0}}{\binom{n+\ell-1}{\ell-1}} & \frac{\binom{\ell-1}{1}}{\binom{n+\ell-2}{\ell-2}} & \frac{\binom{\ell-1}{2}}{\binom{n+\ell-3}{\ell-3}} & \dotsm & \frac{\binom{\ell-1}{\ell-2}}{\binom{n+1}{1}} & \frac{\binom{\ell-1}{\ell-1}}{\binom{n}{0}}\\
\frac{1}{\binom{n+\ell+1}{\ell+1}} & \frac{\binom{\ell}{0}}{\binom{n+\ell}{\ell}} & \frac{\binom{\ell}{1}}{\binom{n+\ell-1}{\ell-1}} & \frac{\binom{\ell}{2}}{\binom{n+\ell-2}{\ell-2}} & \dotsm & \frac{\binom{\ell}{\ell-2}}{\binom{n+2}{2}} & \frac{\binom{\ell}{\ell-1}}{\binom{n+1}{1}}
\end{vmatrix}\\
&=\begin{vmatrix}
A_{\ell+1,1}(n) & B_{\ell+1,\ell}(n)
\end{vmatrix}_{(\ell+1)\times(\ell+1)}
\end{align*}
for $n\in\mathbb{N}$ and $\ell\in\mathbb{N}_0$, in which the matrices $A_{\ell+1,1}(n)$ and $B_{\ell+1,\ell}(n)$ are defined by
\begin{equation*}
A_{\ell+1,1}(n)=\begin{pmatrix}
\alpha_{i,j}(n)
\end{pmatrix}_{1\le i\le \ell+1,j=1}, \quad
\alpha_{i,1}(n)=\frac{1}{\binom{n+i}{i}}
\end{equation*}
and
\begin{equation*}
B_{\ell+1,\ell}(n)=\begin{pmatrix}
\beta_{i,j}(n)
\end{pmatrix}_{1\le i\le \ell+1,1\le j\le \ell}, \quad \beta_{i,j}=\begin{dcases}0, & \ell\ge j>i\in\mathbb{N};\\
\frac{\binom{i-1}{j-1}}{\binom{n+i-j}{i-j}}, & 1\le j\le i\le \ell+1.\end{dcases}
\end{equation*}
\par
To some extent, the determinantal expression~\eqref{Ak=determ} is simpler and more beautiful than the formula~\eqref{A(k)Bell-Eq}.
\item
For $n\in\mathbb{N}$, the ratio
\begin{equation}\label{R(u)(n0)-log}
R_{n,0}(x)=
\begin{dcases}
\frac{\ln T_{n-1}[\te^x]}{x}, & x\ne0\\
\frac{1}{n+1}, & x=0
\end{dcases}
\end{equation}
is increasing on $(-\infty,\infty)$.
\item
The normalized remainder $T_{n-1}[\te^x]$ is increasing and logarithmically convex on $(-\infty,\infty)$. This conclusion generalizes the first result in~\cite[Theorem~2.1]{Guo-Qi-Filomat-2011-May-12.tex}.
\item
For $n\in\mathbb{N}$ and $x\in\mathbb{R}$, we have
\begin{equation}\label{log-exp-last-proof-ineq}
\Biggl[\sum_{j=0}^{\infty}\frac{1}{\binom{n+j+2}{n}}\frac{x^{j}}{j!}\Biggr] \Bigg[\sum_{j=0}^{\infty}\frac{1}{\binom{n+j}{n}}\frac{x^{j}}{j!}\Biggr] \ge\Biggl[\sum_{j=0}^{\infty}\frac{1}{\binom{n+j+1}{n}}\frac{x^{j}}{j!}\Biggr]^2.
\end{equation}
\end{enumerate}

\subsubsection{Ratio and Monotonicity}
In the paper~\cite{exp-norm-tail-ratio.tex}, we gained the following results:
\begin{enumerate}[1)]
\item
For $n\in\mathbb{N}$ and $x\in\mathbb{R}$, the normalized remainder $T_{n-1}[\te^x]$ is positive.
\item
For $n\in\mathbb{N}$, the ratio $\frac{T_{n}[\te^x]}{T_{n-1}[\te^x]}$ is decreasing in $x\in(0,\infty)$.
\item
For $n\in\mathbb{N}$, the normalized remainder $T_{n-1}[\te^x]$ is an absolutely monotonic function in $x\in\mathbb{R}$.
\end{enumerate}

\subsubsection{Inequality and Decreasing Monotonicity of Ratio}
In~\cite[Theorem~1]{exp-tail-ratio.tex}, we verified the strict inequality
\begin{multline}\label{log-exp-last-inverse-ineq}
\Biggl[\sum_{j=0}^{\infty}\frac{1}{\binom{j+n+2}{n+1}}\frac{x^{j}}{j!}\Biggr] \Biggl[\sum_{j=0}^{\infty}\frac{1}{\binom{j+n}{n}}\frac{x^{j}}{j!}\Biggr]\\
<\Biggl[\sum_{j=0}^{\infty}\frac{1}{\binom{j+n+1}{n+1}}\frac{x^{j}}{j!}\Biggr] \Biggl[\sum_{j=0}^{\infty}\frac{1}{\binom{j+n+1}{n}}\frac{x^{j}}{j!}\Biggr]
\end{multline}
for $x\in\mathbb{R}$ and $n\in\mathbb{N}_0$. In~\cite[Theorem~2]{exp-tail-ratio.tex}, by virtue of the inequality~\eqref{log-exp-last-inverse-ineq}, we proved that the ratio $\frac{T_{n+1}[\te^x]}{T_{n}[\te^x]}$ for $n\in\mathbb{N}_0$ is decreasing in $x\in\mathbb{R}$.
\par
Theorem~2 in~\cite{exp-tail-ratio.tex} extends~\cite[Theorem~1]{exp-norm-tail-ratio.tex} from $x\in(0,\infty)$ to $x\in\mathbb{R}$.

\subsubsection{Elegant Proofs}
In the paper~\cite{very-short-exp-nt.tex}, by establishing the integral representation
\begin{equation}\label{Tnexp(x)-int}
T_{n}[\te^x]=(n+1)\int_0^1v^{n}\te^{x(1-v)}\td v
\end{equation}
for $n\in\mathbb{N}_0$ and $x\in\mathbb{R}$, and making use of two monotonicity rules in~\cite[pp.~10--11, Theorem~1.25]{AVV-1997} and in~\cite[Lemma~9]{Alice-y-x-conj-one.tex}, \cite[Lemma~12]{manifold-gamma.tex}, and~\cite[Remark~7.2]{Ouimet-LCM-BKMS.tex}, respectively, we concisely and elegantly proved the following properties of the normalized remainder $T_{n}[\te^x]$ once again:
\begin{enumerate}[1)]
\item
For $n\in\mathbb{N}_0$, the normalized remainder $T_{n}[\te^x]$ is positive in $x\in\mathbb{R}$; see also~\cite[Section~2]{exp-norm-tail-ratio.tex}.
\item
For $n\in\mathbb{N}_0$, the normalized remainder $T_n[\te^x]$ is an increasing and logarithmically convex function of $x\in\mathbb{R}$; see also~\cite[Corollary~1]{log-exp-expan-Sym.tex}.
\item
For $n\in\mathbb{N}_0$, the normalized remainder $T_{n}[\te^x]$ is an absolutely monotonic function in $x\in\mathbb{R}$; see also~\cite[Theorem~2]{exp-norm-tail-ratio.tex}.
\item
For $n\in\mathbb{N}$, the function $R_{n,0}(x)$ defined by~\eqref{R(u)(n0)-log} is increasing in $x\in\mathbb{R}$; see also~\cite[Theorem~2]{log-exp-expan-Sym.tex}.
\item
For $x\in\mathbb{R}$ and $n\in\mathbb{N}_0$, the inequalities~\eqref{log-exp-last-proof-ineq} and~\eqref{log-exp-last-inverse-ineq} with equality are true; see also~\cite[Corollary~2]{log-exp-expan-Sym.tex} and~\cite[Theorem~1]{exp-tail-ratio.tex}.
\item
For $n\in\mathbb{N}_0$, the ratio $\frac{T_{n+1}[\te^x]}{T_{n}[\te^x]}$ is decreasing in $x\in\mathbb{R}$; see also~\cite[Theorem~1]{exp-norm-tail-ratio.tex} and~\cite[Theorem~2]{exp-tail-ratio.tex}.
\end{enumerate}

\subsubsection{More General Results}
In the paper~\cite{ser-int-norm-tail.tex}, by establishing the integral representations
\begin{equation}\label{S(mn)(u)-ser-defn}
S_{m,n}(x)=\sum_{j=0}^{\infty}\frac{1}{\binom{j+m+n}{n}}\frac{x^{j}}{j!}
=n\int_0^1v^{n-1}(1-v)^m\te^{x(1-v)}\td v
\end{equation}
and
\begin{equation*}
S_{\alpha,\beta}(x)=\sum_{j=0}^{\infty}\frac{1}{\binom{j+\alpha+\beta}{\beta}}\frac{x^j}{j!}
=\beta\int_0^1v^{\beta-1}(1-v)^\alpha\te^{u(1-v)}\td v
\end{equation*}
for $m,n\in\mathbb{N}_0$, $\Re(\alpha)>-1$, $\Re(\beta)>0$, and $x\in\mathbb{R}$, 
and by finding the relation
\begin{equation*}
S_{m,n+1}(x)=(T_n[\te^x])^{(m)}
\end{equation*}
for $x\in\mathbb{R}$ and $m,n\in\mathbb{N}_0$, we reformulated the inequalities~\eqref{log-exp-last-proof-ineq}, \eqref{log-exp-last-inverse-ineq}, and
\begin{equation*}
\Biggl[\sum_{j=0}^{\infty}\frac{j+1}{(j+n+2)!}x^{j}\Biggr] \Biggl[\sum_{j=0}^{\infty}\frac{j+1}{(j+n)!}x^{j}\Biggr]
<\Biggl[\sum_{j=0}^{\infty}\frac{j+1}{(j+n+1)!}x^{j}\Biggr]^2
\end{equation*}
for $n\in\mathbb{N}$ and $x\in\mathbb{R}$ in~\cite[Lemma~3]{exp-tail-ratio.tex} as integral forms and in terms of $S_{m,n}(x)$, alternatively proved the absolute monotonicity in $x\in\mathbb{R}$ of the normalized remainder $T_n[\te^x]$ for $n\in\mathbb{N}_0$, and posed the following problems.

\begin{prob}[{\cite[Problem~1]{ser-int-norm-tail.tex}}]\label{2limits-int}
For $m,n\in\mathbb{N}$ and $x\in\mathbb{R}$, calculate the limits
\begin{equation*}
\lim_{n\to\infty}\frac{S_{m,n}(x)} {S_{m-1,n}(x)}
=\lim_{n\to\infty}\frac{\int_0^1v^{n-1}(1-v)^m\te^{-x v}\td v} {\int_0^1v^{n-1}(1-v)^{m-1}\te^{-x v}\td v}
\end{equation*}
and
\begin{equation*}
\lim_{m\to\infty}\frac{S_{m,n}(x)} {S_{m-1,n}(x)}
=\lim_{m\to\infty}\frac{\int_0^1v^{n-1}(1-v)^m\te^{-x v}\td v} {\int_0^1v^{n-1}(1-v)^{m-1}\te^{-x v}\td v}.
\end{equation*}
\end{prob}

\begin{prob}[{\cite[Problem~2]{ser-int-norm-tail.tex}}]
For $\Re(\alpha)>0$, $\Re(\beta)>0$, and $x\in\mathbb{R}$, compute the limits
\begin{equation*}
\lim_{\beta\to\infty}\frac{S_{\alpha,\beta}(x)} {S_{\alpha-1,\beta}(x)}
=\lim_{\beta\to\infty}\frac{\int_0^1v^{\beta-1}(1-v)^\alpha\te^{-x v}\td v} {\int_0^1v^{\beta-1}(1-v)^{\alpha-1}\te^{-x v}\td v}
\end{equation*}
and
\begin{equation*}
\lim_{\alpha\to\infty}\frac{S_{\alpha,\beta}(x)} {S_{\alpha-1,\beta}(x)}
=\lim_{\alpha\to\infty}\frac{\int_0^1v^{\beta-1}(1-v)^\alpha\te^{-x v}\td v} {\int_0^1v^{\beta-1}(1-v)^{\alpha-1}\te^{-x v}\td v}.
\end{equation*}
\end{prob}

\subsection{Normalized Remainders of Inverse Trigonometric Functions}
In~\cite{LiuXL-Arc-three.tex}, we discussed the normalized remainders
\begin{align*}
T_{2n}[\arcsin x]&=T_{2n}[\arccos x]\\
&=
\begin{dcases}
\frac{[(2n)!!]^2}{(2n)!}\frac{2n+1}{x^{2n+1}} \Biggl(\arcsin x-\sum_{k=0}^{n-1}\frac{(2k)!}{[(2k)!!]^2}\frac{x^{2k+1}}{2k+1}\Biggr), & x\ne0;\\
1, & x=0,
\end{dcases}\\
T_{2n}[\arctan x]&=
\begin{dcases}
(-1)^{n}\frac{2n+1}{x^{2n+1}}\Biggl[\arctan x-\sum_{k=0}^{n-1}(-1)^k\frac{x^{2k+1}}{2k+1}\Biggr], & x\ne0;\\
1, & x=0,
\end{dcases}
\end{align*}
and
\begin{equation*}
T_{2n}[\arctanh x]=
\begin{dcases}
\frac{2n+1}{x^{2n+1}}\Biggl(\arctanh x-\sum_{k=0}^{n-1}\frac{x^{2k+1}}{2k+1}\Biggr), & x\ne0;\\
1, & x=0
\end{dcases}
\end{equation*}
for $n\in\mathbb{N}_0$ and $x\in(-1,1)$. Among other findings, we derived the following conclusions:
\begin{enumerate}[1)]
\item
For $n\in\mathbb{N}_0$, the normalized remainder $T_{2n}[\arcsin x]$ is a logarithmically absolutely monotonic function on $(0,1)$ and a logarithmically completely monotonic function on $(-1,0)$.
\item
For $n\in\mathbb{N}_0$, the logarithmic function $\ln T_{2n}[\arcsin x]$ is absolutely monotonic on $(0,1)$ and completely monotonic on $(-1,0)$.
\item
For $n\in\mathbb{N}_0$, the ratio $\frac{T_{2n+2}[\arcsin x]}{T_{2n}[\arcsin x]}$ is an absolutely monotonic function on $(0,1)$ and a completely monotonic function on $(-1,0)$.
\item
For $n\in\mathbb{N}_0$, the normalized remainder $T_{2n}[\arctan x]$ is decreasing in $x\in(0,1)$.
\item
For $n\in\mathbb{N}_0$, the normalized remainder
$$
T_{2n}[\arctan(x\ti)]=T_{2n}[\arctanh x]
$$
is a logarithmically absolutely monotonic function on $(0,1)$ and a logarithmically completely monotonic function on $(-1,0)$.
\item
For $n\in\mathbb{N}_0$, the logarithmic function
$$
\ln T_{2n}[\arctan(x\ti)]=\ln T_{2n}[\arctanh x]
$$
is absolutely monotonic on $(0,1)$ and completely monotonic on $(-1,0)$, where $\ti=\sqrt{-1}\,$ is the imaginary unit.
\item
Let $x_0\in\bigl(-\frac{1}{2},-\frac{2}{5}\bigr)$ be the unique root of the equation
\begin{equation}\label{zero-point}
\sqrt{1-x^2}\,+x\arccos x=0, \quad x\in(-1,1).
\end{equation}
Then the function $\arccos x$ is logarithmically concave on $(-1,x_0)$ and logarithmically convex on $(x_0,1)$.
\item
The logarithmic functions
\begin{equation*}
\ln\biggl(\frac{2}{\pi}\arccos x\biggr), \quad \ln T_{2n}[\arctan x], \quad \ln T_{2n}[\arcsin x]
\end{equation*}
were expanded into the Maclaurin power series.
\end{enumerate}

\subsection{Normalized Remainders of a Function Involved in an Integral Representation of Reciprocal of Gamma Function}
In~\cite[p.~71, Eq.~(3.38)]{Temme-96-book}, we find the integral representation
\begin{equation*}
\frac{1}{\Gamma(z)}=\frac{\te^z z^{1-z}}{\pi}\int_{0}^{\pi}\te^{-z\Phi(\theta)}\td\theta, \quad \Re(z)\ge0,
\end{equation*}
where the classical Euler gamma function $\Gamma(z)$ can be defined~\cite[Chapter~3]{Temme-96-book} by
\begin{equation*}
\Gamma(z)=\lim_{n\to\infty}\frac{n!n^z}{\prod_{k=0}^n(z+k)}, \quad z\in\mathbb{C}\setminus\{0,-1,-2,\dotsc\}
\end{equation*}
and
\begin{equation}
\begin{aligned}\label{Phi-series}
\Phi(\theta)&=1-\frac{\theta}{\tan\theta}+\ln\frac{\theta}{\sin\theta}\\
&=\sum_{j=1}^{\infty}\frac{2j+1}{2j}|B_{2j}|\frac{(2\theta)^{2j}}{(2j)!}\\
&=\frac{1}{2}\theta^2+\frac{1}{36}\theta^4+\frac{1}{405}\theta^6+\frac{1}{4200}\theta^8+\frac{1}{42525}\theta^{10}+\dotsm, \quad |\theta|<\pi.
\end{aligned}
\end{equation}
\par
By the way, combining the series expansion~\eqref{ln-sin-series} with the series expansion~\eqref{Phi-series}, we derive
\begin{equation*}
\frac{\theta}{\tan\theta}=1-\sum_{k=1}^{\infty}|B_{2k}|\frac{(2\theta)^{2k}}{(2k)!}, \quad |\theta|<\pi.
\end{equation*}
This can also be derived by the Laurent series expansion of $\cot\theta$ in the handbooks~\cite[p.~75, Entry~4.3.70]{abram} and~\cite[p.~117, Entry~4.19.6]{NIST-HB-2010}; see also~\cite[Remark~1]{gam-recip-int-phi-sym.tex}. As a result, the function $1-\frac{\theta}{\tan\theta}$ is absolutely monotonic on $(0,\pi)$ and is completely monotonic on $(-\pi,0)$.

\subsubsection{Normalized Remainder of $\Phi(\theta)$}
In the paper~\cite{gam-recip-int-phi-sym.tex}, we considered the function $\Phi(\theta)$ and its normalized remainder
\begin{multline*}
T_{2n+1}[\Phi(\theta)]=\\
\begin{dcases}
\frac{2n+2}{2n+3}\frac{1}{|B_{2n+2}|} \frac{(2n+2)!}{(2\theta)^{2n+2}} \Biggl[\Phi(\theta)-\sum_{j=1}^{n}\frac{2j+1}{2j}|B_{2j}|\frac{(2\theta)^{2j}}{(2j)!}\Biggr], & \theta\ne0\\
1, & \theta=0
\end{dcases}
\end{multline*}
for $n\in\mathbb{N}_0$ and $\theta\in(-\pi,\pi)$. Our main results are as follows:
\begin{enumerate}[1)]
\item
The normalized remainder $T_{2n+1}[\Phi(\theta)]$ for $n\in\mathbb{N}_0$ is positive in $\theta\in(-\pi,\pi)$, decreasing in $\theta\in(-\pi,0)$, and increasing in $\theta\in(0,\pi)$.
\item
The normalized remainder $T_{2n+1}[\Phi(\theta)]$ for $n\in\mathbb{N}_0$ is convex in $\theta\in(-\pi,\pi)$.
\item
The normalized remainder $T_{2n+1}[\Phi(\theta)]$ for $n\in\mathbb{N}_0$ is an absolutely monotonic function in $\theta\in(0,\pi)$ and a completely monotonic function in $\theta\in(-\pi,0)$.
\item
The normalized remainder $T_{2n+1}[\Phi(\theta)]$ for $n\in\mathbb{N}_0$ is logarithmically convex in $\theta\in(-\pi,\pi)$.
\item
The ratio $\frac{T_{2n+3}[\Phi(\theta)]}{T_{2n+1}[\Phi(\theta)]}$ for $n\in\mathbb{N}_0$ is decreasing in $\theta\in(-\pi,0)$ and increasing in $\theta\in(0,\pi)$.
\end{enumerate}

\subsubsection{Stronger Properties of Normalized Remainder of $\Phi(\theta)$}
In the paper~\cite{JMI-5203.tex}, we generalized the main results in~\cite{gam-recip-int-phi-sym.tex} and presented the following conclusions:
\begin{enumerate}[1)]
\item
For $n\in\mathbb{N}_0$, the normalized remainder $T_{2n+1}[\Phi(\theta)]$ is a logarithmically absolutely monotonic function in $\theta\in(0,\pi)$ and a logarithmically completely monotonic function in $\theta\in(-\pi,0)$.
\item
For $n\in\mathbb{N}_0$, the function $\frac{1}{\theta}\frac{\td\,\ln T_{2n+1}[\Phi(\theta)]}{\td \theta}$ is an absolutely monotonic function in $(0,\pi)$ and a completely monotonic function in $(-\pi,0)$.
\item
For $n\in\mathbb{N}_0$, the normalized remainder $T_{2n+1}[\Phi(\theta)]$ and $\frac{1}{\theta}\frac{\td\,\ln T_{2n+1}[\Phi(\theta)]}{\td \theta}$ can be extended analytically into the complex $z$-plane and are analytic in the disc $|z|<\pi$.
\item
For $n\in\mathbb{N}_0$, the ratio $\frac{T_{2n+3}[\Phi(\theta)]}{T_{2n+1}[\Phi(\theta)]}$
\begin{enumerate}[1)]
\item
is an absolutely monotonic function in $\theta\in(0,\pi)$ and a completely monotonic function in $\theta\in(-\pi,0)$,
\item
can be extended analytically into the complex $z$-plane and is analytic in the disc $|z|<\pi$.
\end{enumerate}
\item
The function $\frac{1}{\Phi(\theta)}$ for $0<|\theta|<\pi$ was expanded into a Laurent series.
\end{enumerate}
\par
The absolute and complete monotonicity of the ratio $\frac{T_{2n+3}[\Phi(\theta)]}{T_{2n+1}[\Phi(\theta)]}$ discovered in~\cite[Theorem~2]{JMI-5203.tex} and the Laurent series of the function $\frac{1}{\Phi(\theta)}$ for $0<|\theta|<\pi$ expanded in~\cite[Theorem~3]{JMI-5203.tex} imply that the function $\frac{2}{\theta^2}-\frac{1}{\Phi(\theta)}$ is absolutely monotonic on $(0,\pi)$ and completely monotonic on $(-\pi,0)$.

\subsection{Normalized Remainders of Complete Elliptic Integrals}
For complex numbers $a, b, c\in\mathbb{C}$ with $c\ne 0, -1,-2,\dotsc$, the Gauss hypergeometric function ${}_{2}F_1$ is defined~\cite{Gauss-Milovanovic-Qi.tex, Gauss-hyperg-Int.tex} by
\begin{equation*}
{}_{2}F_1(a,b;c;z)=\sum_{n=0}^\infty\frac{(a)_n(b)_n}{(c)_n}\frac{z^n}{n!},\quad |z|<1,
\end{equation*}
where $(z)_n$ for $z\in\mathbb{C}$ is the Pochhammer symbol, also known as the rising factorial and the shifted factorial, defined by
\begin{equation*}
(z)_n=\prod_{k=0}^{n-1}(z+k)=
\begin{cases}
z(z+1)\dotsm(z+n-1), & n\in\mathbb{N};\\
1, & n=0.
\end{cases}
\end{equation*}
\par
For given $x\in[0,1]$, the complete elliptic integral of the second kind $\E(x)$ are defined~\cite[p.~315]{Temme-96-book} by
\begin{equation*}
\E(x)=\int_{0}^{\pi/2}\sqrt{1-x^2\sin^2\theta}\,\td\theta
=\int_{0}^{1}\frac{\sqrt{1-x^2t^2}\,}{\sqrt{1-t^2}\,}\td t
\end{equation*}
and can be represented~\cite[p.~128]{Temme-96-book} by
\begin{equation*}
\begin{dcases}
\E(x)=\frac{\pi}{2}{\,}_{2}F_1\biggl(-\frac{1}{2},\frac{1}{2};1;x^2\biggr);\\
\E(0)=\frac{\pi}{2};\\
\E(1)=1.
\end{dcases}
\end{equation*}
\par
Let
\begin{equation*}
\mathcal{H}(x)=\frac{1-{\,}_{2}F_1\bigl(-\frac12,\frac12;1;x\bigr)}{1-(1-x){\,}_{2}F_1\bigl(\frac12,1;\frac32;x\bigr)}
=\frac{1-\frac{2}{\pi}\E(\sqrt{x}\,)}{1+(x-1) \frac{\arctanh\sqrt{x}\,}{\sqrt{x}\,}}
\end{equation*}
and
\begin{align*}
\mathcal{G}(x)&=\frac{\frac{\pi}{2}\bigl[{\,}_{2}F_1\bigl(-\frac12,\frac12;1;x\bigr)-1\bigr] +x{\,}_{2}F_1\bigl(\frac12,1;\frac32;x\bigr)}{x{\,}_{2}F_1(1,1;2;x)}\\
&=\frac{\frac\pi2-\sqrt{x}\, \arctanh\sqrt{x}\,-\E(\sqrt{x}\,)}{\ln (1-x)}.
\end{align*}

\subsubsection{First Proof of Absolute Monotonicity}
In the paper~\cite[Remarks~2 and~4]{JMI-5046.tex}, we proved that the normalized remainders $T_n[\mathcal{H}(x)]$ and $T_n[\mathcal{G}(x)]$ for $n\in\mathbb{N}_0$ are absolutely monotonic on $(0,1)$.
\subsubsection{Second Proof of Absolute Monotonicity}
In the paper~\cite[Theorems~3 and~4]{ratio-tan-squ.tex}, we employed a new method to concisely prove the absolute monotonicity of the normalized remainders $T_n[\mathcal{H}(x)]$ and $T_n[\mathcal{G}(x)]$ for $n\in\mathbb{N}_0$ on $(0,1)$ once again.

\subsection{Normalized Remainder of Power of Inverse Sine Function}
In~\cite[Theorem~1]{AADM-3164.tex}, \cite[Theorem~2.1]{AIMS-Math20210491.tex}, \cite[Section~4]{series-arccos-up-42.tex}, and~\cite[Section~6]{dema-D-22-00034.tex}, we discovered the Maclaurin power series expansion
\begin{equation}\label{arcsin-series-expansion-unify}
\biggl(\frac{\arcsin x}{x}\biggr)^n
=1+\sum_{\ell=1}^{\infty} \frac{\mathcal{Q}(n,2\ell)}{\binom{n+2\ell}{2\ell}} \frac{(2x)^{2\ell}}{(2\ell)!}, \quad |x|<1
\end{equation}
for $n\in\mathbb{N}$ several times by different ideas and techniques, where
\begin{equation}\label{Q(m-n)-sum-dfn}
\mathcal{Q}(n,2\ell)=(-1)^\ell\sum_{k=0}^{2\ell} \binom{n+k-1}{n-1} s(n+2\ell-1,n+k-1)\biggl(\frac{n+2\ell-2}{2}\biggr)^{k}
\end{equation}
for $n,\ell\in\mathbb{N}$ and the Stirling numbers of the first kind $s(\ell+n,n)$ for $\ell,n\in\mathbb{N}_0$ can be analytically generated~\cite[Theorem~3.14]{Mansour-Schork-B2016} by
\begin{equation}\label{Stirl-No-First-GF}
\biggl[\frac{\ln(1+x)}{x}\biggr]^n=\sum_{\ell=0}^\infty \frac{s(\ell+n,n)}{\binom{\ell+n}{n}}\frac{x^{\ell}}{\ell!}, \quad |x|<1.
\end{equation}
\par
The power series expansion~\eqref{arcsin-series-expansion-unify} is the unification of the Maclaurin power series
\begin{align*}
\frac{\arcsin x}{x}&=\sum_{\ell=0}^{\infty}\frac{[(2\ell-1)!!]^2}{\binom{2\ell+1}{1}}\frac{x^{2\ell}}{(2\ell)!},\\
\biggl(\frac{\arcsin x}{x}\biggr)^2&=\sum_{\ell=0}^{\infty} \frac{[(2\ell)!!]^2}{\binom{2\ell+2}{2}} \frac{x^{2\ell}}{(2\ell)!},\\
\biggl(\frac{\arcsin x}{x}\biggr)^3&=\sum_{\ell=0}^{\infty}\frac{[(2\ell+1)!!]^2}{\binom{2\ell+3}{3}} \Biggl[\sum_{k=0}^{\ell}\frac{1}{(2k+1)^2}\Biggr]\frac{x^{2\ell}}{(2\ell)!},\\
\biggl(\frac{\arcsin x}{x}\biggr)^4&=\sum_{\ell=0}^{\infty}\frac{[(2\ell+2)!!]^2}{\binom{2\ell+4}{4}}\Biggl[\sum_{k=0}^{\ell}\frac{1}{(2k+2)^2}\Biggr] \frac{x^{2\ell}}{(2\ell)!},\\
\biggl(\frac{\arcsin x}{x}\biggr)^5&=\frac{1}{2}\sum_{\ell=0}^{\infty}\frac{[(2\ell+3)!!]^2}{\binom{2\ell+5}{5}} \Biggl[\Biggl(\sum_{k=0}^{\ell+1}\frac{1}{(2k+1)^2}\Biggr)^2 -\sum_{k=0}^{\ell+1}\frac{1}{(2k+1)^4}\Biggr] \frac{x^{2\ell}}{(2\ell)!}
\end{align*}
for $|x|<1$. These series expansions can be found in the literature~\cite[p.~444]{Apostol-Calculus-v1-1967}, \cite[Remark~6]{AADM-3164.tex}, \cite[Remark~5.2]{AIMS-Math20210491.tex}, \cite[Section~1.1]{arcsine-liu-qi-mdpi.tex}, \cite[Section~1]{Qi-Chen-Lim-RNA.tex}, \cite[Section~1]{arcsin-power-wei.tex}, \cite[pp.~279 and~331--333]{Almost-Impossible-B-2019}, \cite[p.~229]{Zorich-Math-Anal-2015}, and~\cite{Bailey2Moll-Experim2007, Borwein-Chamberland-IJMMS-2007, Spivey-art-2019}, for example.
\par
In~\cite[Theorem~7]{series-arccos-up-42.tex}, see also~\cite[Example~3.1]{integer2Bell2real.tex}, the Maclaurin series expansion~\eqref{arcsin-series-expansion-unify} was generalized as
\begin{equation*}
\biggl(\frac{\arcsin x}{x}\biggr)^r
=1+\sum_{\ell=1}^{\infty} \Biggl[\sum_{k=1}^{2\ell} \frac{(-r)_k}{(2\ell+k)!} \sum_{j=1}^{k}(-1)^{j}\binom{2\ell+k}{k-j}\mathcal{Q}(j,2\ell)\Biggr] (2x)^{2\ell}
\end{equation*}
for $r\in\mathbb{R}$ and $|x|<1$.
\par
In the paper~\cite{arcsine-liu-qi-mdpi.tex}, we mainly considered the normalized remainder \small
\begin{multline*}
T_{2n+1}\biggl[\biggl(\frac{\arcsin x}{x}\biggr)^q\biggr]\\
=
\begin{dcases}
\frac{\binom{q+2n+2}{2n+2}}{\mathcal{Q}(q,2n+2)} \frac{(2n+2)!}{(2x)^{2n+2}} \Biggl[\biggl(\frac{\arcsin x}{x}\biggr)^q-1 -\sum_{\ell=1}^n\frac{\mathcal{Q}(q,2\ell)}{\binom{q+2\ell}{2\ell}}\frac{(2x)^{2\ell}}{(2\ell)!}\Biggr], & x\ne0\\
1, & x=0
\end{dcases}
\end{multline*} \normalsize
for $q\in\mathbb{N}$, $n\in\mathbb{N}_0$, and $|x|<1$, where $\mathcal{Q}(q,2\ell)$ is defined by~\eqref{Q(m-n)-sum-dfn}, and acquired the following conclusions:
\begin{enumerate}[1)]
\item
The normalized remainder $T_{2n+1}\bigl[\bigl(\frac{\arcsin x}{x}\bigr)^q\bigr]$ is an even function in $x\in(-1,1)$.
\item
For $n,\ell\in\mathbb{N}$, the central factorial numbers of the first kind $t(n,\ell)$ satisfy
\begin{equation*}
t(n+2\ell-1,n)=0
\end{equation*}
and
\begin{equation*}
\mathcal{Q}(n,2\ell)=(-1)^{\ell}t(n+2\ell,n),
\end{equation*}
where the central factorial numbers of the first kind $t(n,\ell)$ are defined~\cite[pp.~212--217, Section~6.5]{Riordan-B-1968} by
\begin{equation*}
x\prod_{\ell=1}^{n-1}\biggl(x+\frac{n}{2}-\ell\biggr)=\sum_{\ell=1}^{n}t(n,\ell)x^\ell, \quad n\in\mathbb{N}
\end{equation*}
and can be generated~\cite[p.~44]{Liu-GD-Debrecen-2011} by
\begin{equation*}
\Biggl[2\ln\Biggl(\frac{x}{2}+\sqrt{1+\biggl(\frac{x}{2}\biggr)^2}\,\Biggr)\Biggr]^n
=n!\sum_{\ell=n}^{\infty}t(\ell,n)\frac{x^\ell}{\ell!}, \quad n\in\mathbb{N}.
\end{equation*}
Consequently, for $n\in\mathbb{N}$, we have
\begin{equation*}
\biggl(\frac{\arcsin x}{x}\biggr)^n
=1+\sum_{\ell=1}^{\infty}(-1)^{\ell} \frac{t(n+2\ell,n)}{\binom{n+2\ell}{2\ell}} \frac{(2x)^{2\ell}}{(2\ell)!}, \quad |x|<1.
\end{equation*}
\item
For $n\ge\ell\in\mathbb{N}$, the central factorial numbers of the first kind $t(n,\ell)$ can be explicitly computed by
\begin{equation*}
t(n,\ell)=Q(\ell,n-\ell),
\end{equation*}
where $Q(n,k)$ is defined by
\begin{equation*}
Q(n,k)=\sum_{\ell=0}^{k} \binom{n+\ell-1}{n-1} s(n+k-1,n+\ell-1)\biggl(\frac{n+k-2}{2}\biggr)^{\ell}
\end{equation*}
for $n\in\mathbb{N}$ and $k\in\mathbb{N}_0$, with assumptions $Q(1,1)=0$ and $Q(2,0)=1$, satisfying $Q(n,2\ell)=(-1)^\ell\mathcal{Q}(n,2\ell)$.
\item
The central factorial numbers of the first kind $t(2n,2\ell)$ can be explicitly computed by
\begin{equation}\label{t(2n-2ell)-explicit}
t(2n,2\ell)=\sum_{k=2\ell-1}^{2n-1}\binom{k}{2\ell-1}s(2n-1,k)(n-1)^{k-2\ell+1},\quad n\ge\ell\in\mathbb{N}
\end{equation}
and the central factorial numbers of the first kind $t(2n-1,2\ell-1)$ can be explicitly computed by
\begin{equation*}
t(2n+1,2\ell+1)
=\sum_{k=2\ell}^{2n}\binom{k}{2\ell} s(2n,k) \biggl(n-\frac12\biggr)^{k-2\ell}, \quad 0\le\ell\le n,
\end{equation*}
where the factor $0^0$ when $n=1$ in~\eqref{t(2n-2ell)-explicit} is assumed to be $1$ and $s(n,k)$ denotes the Stirling numbers of the first kind generated in~\eqref{Stirl-No-First-GF}.
\item
The central factorial numbers of the first kind $t(n,m)$ can be explicitly computed by
\begin{multline*}
t(2n,2m)=(-1)^{n-1}\Biggl([(n-1)!]^2\sum_{j=m-1}^{n-1}\frac{|s(j,m-1)|}{j!}\\
+\sum_{i=1}^{n-1}(-1)^{i} \Biggl[\prod_{\ell=1}^{n-1}\bigl(i+\ell^2\bigr)\Biggr] \Biggl[\sum_{j=m-1}^{n-1} \binom{j}{i}\frac{|s(j,m-1)|}{j!}\Biggr]\Biggr)
\end{multline*}
for $n\ge m\in\mathbb{N}$ and
\begin{multline*}
t(2n+1,2m+1)
=\frac{(-1)^{n}}{4^{n-m}}\Biggl([(2n-1)!!]^2\sum_{j=m}^{n}\frac{|s(j,m)|}{j!}\\
+\sum_{i=1}^{n}(-1)^{i} \Biggl[\prod_{\ell=1}^{n}\bigl[i+(2\ell-1)^2\bigr]\Biggr] \Biggl[\sum_{j=m}^{n} \binom{j}{i}\frac{|s(j,m)|}{j!}\Biggr]\Biggr)
\end{multline*}
for $n\ge m\in\mathbb{N}_0$.
\item
For $n\ge\ell\in\mathbb{N}$, the central factorial numbers of the first kind $t(n,\ell)$ can be computed by
\begin{align*}
t(2n,2\ell)&=(-1)^{n-\ell}\sum_{1\le j_1<j_2<\dotsm<j_{n-\ell}\le n-1} \prod_{m=1}^{n-\ell}j_{m}^2,\\
t(2n-1,2\ell-1)&=(-1)^{n-\ell}\sum_{1\le j_1<j_2<\dotsm<j_{n-\ell}\le n-1} \prod_{m=1}^{n-\ell}\biggl(j_{m}-\frac{1}{2}\biggr)^2,
\end{align*}
and
\begin{equation*}
t(n,\ell)=(-1)^{n-\ell}\sum_{1\le j_1<j_2<\dotsm<j_{n-\ell}\le n-1} \prod_{k=1}^{n-\ell}\biggl(j_k-\frac{n}{2}\biggr).
\end{equation*}
Consequently, for $n\ge \ell\in\mathbb{N}$, we have
\begin{equation*}
(-1)^{n-\ell}t(2n,2\ell)>0 \quad\text{and}\quad (-1)^{n-\ell}t(2n-1,2\ell-1)>0.
\end{equation*}
\item
For $n\in\mathbb{N}$, let the integer Vandermonde matrix
\begin{equation*}
V(n)=
\begin{pmatrix}
 1 & 0 & 0 & \dotsm & 0 & 0\\
 1 & 1 & 1 & \dotsm & 1 & 1\\
 1 & 2 & 2^{2} & \dotsm & 2^{n-1} & 2^{n}\\
 \vdots & \vdots & \vdots & \ddots & \vdots & \vdots\\
 1 & n-1 & (n-1)^{2} & \dotsm & (n-1)^{n-1} & (n-1)^{n}\\
 1 & n & n^{2} & \dotsm & n^{n-1} & n^{n}
\end{pmatrix}_{(n+1)\times(n+1)}.
\end{equation*}
Then its inverse matrix is
\begin{equation*}
V^{-1}(n)=
\begin{pmatrix}
1 & \boldsymbol{0}(n)\\
-U(n)\boldsymbol{1}^{\textup{T}}(n) & U(n)
\end{pmatrix},
\end{equation*}
where $\textup{T}$ stands for the transpose notation, and
\begin{gather*}
\boldsymbol{0}(n)=(\overbrace{0,0,\dotsc,0}^{n}), \quad \boldsymbol{1}(n)=(\overbrace{1,1,\dotsc,1}^{n}),\\
U(n)=
\begin{pmatrix}
u_{k,i}(n)
\end{pmatrix}_{1\le k\le n, 1\le i\le n}
=
\begin{pmatrix}\displaystyle
\frac{1}{i!}\sum_{j=k}^{n}\frac{(-1)^{j-i}}{(j-i)!}s(j,k)
\end{pmatrix}_{1\le k\le n, 1\le i\le n}.
\end{gather*}
\item
For $n,\ell\in\mathbb{N}$, the quantity $\mathcal{Q}(n,2\ell)$ defined by~\eqref{Q(m-n)-sum-dfn} is positive.
\item
For $1\le \ell\le n\in\mathbb{N}$, we have
\begin{equation*}
(-1)^\ell \Biggl((n!)^2\sum_{j=\ell}^{n}\frac{|s(j,\ell)|}{j!} +\sum_{i=1}^{n} (-1)^{i}\Biggl[\sum_{j=\ell}^{n}\binom{j}{i}\frac{|s(j,\ell)|}{j!}\Biggr] \prod_{k=1}^{n}\bigl(i+k^2\bigr)\Biggr)>0
\end{equation*}
and
\begin{multline*}
(-1)^\ell\Biggl([(2n-1)!!]^2\sum_{j=\ell}^{n}\frac{|s(j,\ell)|}{j!} \\
+\sum_{i=1}^{n} (-1)^{i}\Biggl[\sum_{j=\ell}^{n}\binom{j}{i}\frac{|s(j,\ell)|}{j!}\Biggr] \prod_{k=1}^{n}\bigl[i+(2k-1)^2\bigr]\Biggr)>0.
\end{multline*}
\item
For $1\le q\le5$ and $n\in\mathbb{N}_0$, the normalized remainder $T_{2n+1}\bigl[\bigl(\frac{\arcsin x}{x}\bigr)^q\bigr]$ is positive, increasing, convex, absolutely monotonic, and logarithmically convex in $x\in[0,1)$.
\item
The logarithm $\ln T_{2n+1}\bigl[\bigl(\frac{\arcsin x}{x}\bigr)^q\bigr]$ was expanded into a Maclaurin series.
\item
The partial Bell polynomials, also known as the Bell polynomials of the second kind, are denoted and defined in~\cite[Definition~11.2]{Charalambides-book-2002} and~\cite[p.~134, Theorem~A]{Comtet-Combinatorics-74} by
\begin{multline}\label{Bell2nd-Dfn-Eq}
B_{n,k}(x_1,x_2,\dotsc,x_{n-k+1})\\
=\sum_{\substack{1\le j\le n-k+1,\,\ell_j\in\{0\}\cup\mathbb{N}\\ \sum_{j=1}^{n-k+1}j\ell_j=n,\,
\sum_{j=1}^{n-k+1}\ell_j=k}}\frac{n!}{\prod_{j=1}^{n-k+1}\ell_j!} \prod_{j=1}^{n-k+1}\biggl(\frac{x_j}{j!}\biggr)^{\ell_j}.
\end{multline}
Then the identities
\begin{align*}
\frac{\bell_{2n+3,1}\bigl(1, 0, 1, 0, 9, 0,\dotsc,[(2n-1)!!]^2, 0, [(2n+1)!!]^2\bigr)} {\bell_{2n+1,1}\bigl(1, 0, 1, 0, 9, 0,\dotsc,[(2n-3)!!]^2, 0, [(2n-1)!!]^2\bigr)}
&=(1+2n)^2,\\
\frac{\bell_{2n+4,2}\bigl(1, 0, 1, 0, 9, 0,\dotsc,[(2n-1)!!]^2, 0, [(2n+1)!!]^2\bigr)} {\bell_{2n+2,2}\bigl(1, 0, 1, 0, 9, 0,\dotsc,[(2n-3)!!]^2, 0, [(2n-1)!!]^2\bigr)}
&=(2+2n)^2,
\end{align*}
\begin{multline*}
\frac{\bell_{2n+5,3}\bigl(1, 0, 1, 0, 9, 0,\dotsc,[(2n-1)!!]^2, 0, [(2n+1)!!]^2\bigr)} {\bell_{2n+3,3}\bigl(1, 0, 1, 0, 9, 0,\dotsc,[(2n-3)!!]^2, 0, [(2n-1)!!]^2\bigr)}\\
=(3+2n)^2+\frac{1}{\sum_{k=0}^{n}\frac{1}{(2k+1)^2}},
\end{multline*}
\begin{multline*}
\frac{\bell_{2n+6,4}\bigl(1, 0, 1, 0, 9, 0,\dotsc,[(2n-1)!!]^2, 0, [(2n+1)!!]^2\bigr)} {\bell_{2n+4,4}\bigl(1, 0, 1, 0, 9, 0,\dotsc,[(2n-3)!!]^2, 0, [(2n-1)!!]^2\bigr)}\\
=(4+2n)^2+\frac{1}{\sum_{k=0}^{n}\frac{1}{(2k+2)^2}},
\end{multline*}
and
\begin{multline*}
\frac{\bell_{2n+7,5}\bigl(1, 0, 1, 0, 9, 0,\dotsc,[(2n-1)!!]^2, 0, [(2n+1)!!]^2\bigr)} {\bell_{2n+5,5}\bigl(1, 0, 1, 0, 9, 0,\dotsc,[(2n-3)!!]^2, 0, [(2n-1)!!]^2\bigr)}\\
=(5+2n)^2 \frac{\bigl[\sum_{k=0}^{n+2}\frac{1}{(2k+1)^2}\bigr]^2 -\sum_{k=0}^{n+2}\frac{1}{(2k+1)^4}} {\bigl[\sum_{k=0}^{n+1}\frac{1}{(2k+1)^2}\bigr]^2 -\sum_{k=0}^{n+1}\frac{1}{(2k+1)^4}}
\end{multline*}
are valid for $n\in\mathbb{N}$.
\item
For $n\in\mathbb{N}$ and $1\le\ell\le5$, the central factorial numbers $t(2n+\ell,\ell)$ satisfy the identities
\begin{align*}
\frac{t(2n+3,1)}{t(2n+1,1)}&=-\frac{(2n+1)^2}{4}, \\
\frac{t(2n+5,3)}{t(2n+3,3)}&=-\frac{1}{4}\Biggl[(2n+3)^2+\frac{1}{\sum_{k=0}^{n}\frac{1}{(2k+1)^2}}\Biggr],\\
\frac{t(2n+4,2)}{t(2n+2,2)}&=-\frac{(2n+2)^2}{4}, \\
\frac{t(2n+6,4)}{t(2n+4,4)}&=-\frac{1}{4}\Biggl[(2n+4)^2+\frac{1}{\sum_{k=0}^{n}\frac{1}{(2k+2)^2}}\Biggr],
\end{align*}
and
\begin{equation*}
\frac{t(2n+7,5)}{t(2n+5,5)}
=-\frac{(2n+5)^2}{4} \frac{\bigl[\sum_{k=0}^{n+2}\frac{1}{(2k+1)^2}\bigr]^2 -\sum_{k=0}^{n+2}\frac{1}{(2k+1)^4}} {\bigl[\sum_{k=0}^{n+1}\frac{1}{(2k+1)^2}\bigr]^2 -\sum_{k=0}^{n+1}\frac{1}{(2k+1)^4}}.
\end{equation*}
\item
For $m,n\in\mathbb{N}$, we have the relation
\begin{multline*}
\frac{t(n+2m+2,n)}{t(n+2m,n)}\\
=-\frac{1}{4} \frac{\bell_{n+2m+2,n}\bigl(1, 0, 1, 0, 9, 0,\dotsc,[(2m-1)!!]^2, 0, [(2m+1)!!]^2\bigr)} {\bell_{n+2m,n}\bigl(1, 0, 1, 0, 9, 0,\dotsc,[(2m-3)!!]^2, 0, [(2m-1)!!]^2\bigr)}.
\end{multline*}
\end{enumerate}
\par
From the Maclaurin power series expansion~\eqref{arcsin-series-expansion-unify} and the result that the quantity $\mathcal{Q}(n,2\ell)$ for $n,\ell\in\mathbb{N}$ is positive, we conclude that the normalized remainder $T_{2n+1}\bigl[\bigl(\frac{\arcsin x}{x}\bigr)^q\bigr]$ for $q\in\mathbb{N}$ and $n\in\mathbb{N}_0$ is an absolutely monotonic function in $x\in[0,1)$ and a completely monotonic function in $x\in(-1,0]$.

\section{Historic Backgrounds in Combinatorics and Number Theory}
In this section, we collected some concepts, which can be reformulated in terms of normalized remainders, in combinatorics and number theory.

\subsection{Generating Functions of Carlitz--Howard's Generalizations of Bernoulli Numbers and Polynomials}\label{sec4.1}
The equations~\eqref{Bernoulli-GF} and~\eqref{Bernoulli-Polyn-GF} demonstrate that the Bernoulli numbers $B_n$ can be generated by the reciprocal $\frac1{T_0[\te^x]}$ and the Bernoulli polynomials $B_n(t)$ can be generated by the reciprocal $\frac1{T_0[(\te^x-1)\te^{-t x}]}$ or $\frac{\te^{x t}}{T_0[\te^x]}$.
\par
The generalized Bernoulli polynomials $B_n^{(\sigma)}(t)$ can be analytically generated~\cite[p.~4]{Temme-96-book} by
\begin{equation*}
\te^{x t}\biggl(\frac{x}{\te^{x}-1}\biggr)^{\sigma}=\sum_{n=0}^{\infty}B^{(\sigma)}_{n}(t)\frac{x^{n}}{n!},\quad |z|<2\pi.
\end{equation*}
The generating function $\te^{x t}\bigl(\frac{x}{\te^{x}-1}\bigr)^{\sigma}$ can be regarded as $\frac{\te^{x t}}{(T_0[\te^x])^\sigma}$.
\par
In the thesis~\cite{Howard-Duke-1966-thesis} and the papers~\cite{Carlitz-MathMag1960, Howard-Duke-1967-599, Howard-MC-80-977, Howard-MC-1977, Howard-Duke-1967-701}, among other things, Carlitz and Howard defined the sequences and polynomials $A_k$, $A_k(t)$, $A_{n,k}$, and $A_{n,k}(t)$ by~\eqref{Howard-numbers-def} and
\begin{align}
\frac{x^2}{2}\frac{\te^{x t}}{\te^x-x-1}&=\sum_{k=0}^{\infty}A_k(t)\frac{x^k}{k!},\label{Howard-Gen-No-T}\\
\frac{x^n}{n!}\frac1{\te^x-\sum_{k=0}^{n-1}\frac{x^k}{k!}}
&=\sum_{k=0}^{\infty}A_{n,k}\frac{x^k}{k!}, \label{Howard-Gen-No-DFN}\\
\label{Howard-Gen-polyn-DFN}
\frac{x^n}{n!}\frac{\te^{x t}}{\te^x-\sum_{k=0}^{n-1}\frac{x^k}{k!}}
&=\sum_{k=0}^{\infty}A_{n,k}(t)\frac{x^k}{k!}
\end{align}
for $n\in\mathbb{N}$, and they examined a lot of algebraic properties of these sequences and polynomials.
It is easy to see that $A_k(0)=A_k$, $A_{1,k}=B_k$, $A_{2,k}=A_k$, $A_{n,k}(0)=A_{n,k}$, $A_{1,k}(0)=B_k$, $A_{2,k}(0)=A_k$, and $A_{1,k}(t)=B_k(t)$ for $k\in\mathbb{N}_0$. The sequence $A_k$ for $k\in\mathbb{N}_0$ was recently investigated in~\cite[Remark~2]{log-exp-expan-Sym.tex} and~\cite[Section~2]{Bell-value-elem-funct.tex}.
In terms of the notation $T_{n-1}[\te^x]$, the equation~\eqref{Howard-Gen-polyn-DFN} can be reformulated as
\begin{equation*}
\frac{\te^{x t}}{T_{n-1}[\te^x]}
=\sum_{k=0}^{\infty}A_{n,k}(t)\frac{x^k}{k!}.
\end{equation*}
Meanwhile, the generating functions of the sequences and polynomials $B_k$, $B_k(t)$, $A_k$, $A_k(t)$, and $A_{n,k}$ for $n\in\mathbb{N}$ in~\eqref{Bernoulli-GF}, \eqref{Bernoulli-Polyn-GF}, \eqref{Howard-numbers-def}, \eqref{Howard-Gen-No-T}, and~\eqref{Howard-Gen-No-DFN} can be respectively reformulated as
\begin{equation*}
\frac{1}{T_0[\te^x]},\quad \frac{\te^{x t}}{T_0[\te^x]},\quad \frac{1}{T_1[\te^x]},\quad \frac{\te^{x t}}{T_1[\te^x]}, \quad \frac{1}{T_{n-1}[\te^x]}.
\end{equation*}
See also~\cite[Remark~4]{Spivey-JIMA.tex}, \cite[Remark~4]{Spivey-ID-AP-JIMA.tex}, \cite[pp.~1036--1037]{very-short-exp-nt.tex}, and~\cite[Remark~2]{JMI-5046.tex}.

\subsection{Generating Functions of Howard's Generalization of Stirling Numbers of the Second Kind}
The equation~\eqref{2Stirl-funct-rew} shows that the Stirling numbers of the second kind $S(k+\ell,\ell)$ for $k,\ell\in\mathbb{N}_0$ can be generated by $(T_0[\te^x])^\ell$. The equation~\eqref{Howard-numbers-def} shows that the Howard numbers $A_k$ can be generated by the reciprocal $\frac{1}{T_1[\te^x]}$.
\par
The $r$-associate Stirling numbers of the second kind $S_r(k,\ell)$ for $k\ge \ell\ge0$ and $r\in\mathbb{N}_0$ are defined~\cite[p.~303, Eq.~(1.2)]{Howard-Fib-1980i4} by
\begin{equation}\label{r-Associate-Stirling}
\Biggl(\te^x-\sum_{k=0}^{r}\frac{x^k}{k!}\Biggr)^\ell=\Biggl(\sum_{k=r+1}^{\infty}\frac{x^k}{k!}\Biggr)^\ell
=\ell!\sum_{k=(r+1)\ell}^{\infty}S_r(k,\ell)\frac{x^k}{k!};
\end{equation}
see also~\cite[Section~1.12]{Bell-value-elem-funct.tex}.
The equation~\eqref{r-Associate-Stirling} can be reformulated as
\begin{equation}\label{T-Associate-Stirling}
\begin{split}
(T_r[\te^x])^\ell&=\Biggl[\frac{(r+1)!}{x^{r+1}}\Biggl(\te^x-\sum_{k=0}^{r}\frac{x^k}{k!}\Biggr)\Biggr]^\ell\\
&=\frac{\ell![(r+1)!]^\ell}{[(r+1)\ell]!} \sum_{k=0}^{\infty}\frac{S_r(k+(r+1)\ell,\ell)}{\binom{k+(r+1)\ell}{k}} \frac{x^{k}}{k!}
\end{split}
\end{equation}
for $\ell,r\in\mathbb{N}_0$. The equation~\eqref{2Stirl-funct-rew} is a special case $r=0$ of~\eqref{T-Associate-Stirling}. Consequently, the integer powers $(T_r[\te^x])^\ell$ for $\ell,r\in\mathbb{N}_0$ can be regarded as the generating functions of the Stirling numbers of the second kind $S(k,\ell)$ and the $r$-associate Stirling numbers of the second kind $S_r(k,\ell)$ for $k\ge \ell\in\mathbb{N}_0$ and $r\in\mathbb{N}_0$; see also~\cite[Remark~4]{Spivey-JIMA.tex}, \cite[Remark~4]{Spivey-ID-AP-JIMA.tex}, and~\cite[Remark~2]{JMI-5046.tex}.
\par
The Stirling polynomials $S_k(x)$ for $k\in\mathbb{N}_0$ are generated~\cite{MR462961, Stirling-Polyn.tex} by
\begin{equation*}
\biggl(\frac{t}{1-\te^{-t}}\biggr)^{x+1}=\sum_{k=0}^\infty S_k(x)\frac{t^k}{k!}.
\end{equation*}
It is clear that
\begin{equation*}
S_k(x)=B_k^{(x+1)}(x+1).
\end{equation*}
The generating function of $S_k(x)$ can be rewritten as
\begin{equation*}
\biggl(\frac{t}{1-\te^{-t}}\biggr)^{x+1}=\biggl(\frac{1}{T_0[\te^{-t}]}\biggr)^{x+1}.
\end{equation*}

\subsection{Generating Functions of Howard's Generalization of Stirling Numbers of the First Kind}
The Stirling numbers of the first kind $s(n,k)$ for $n,k\in\mathbb{N}_0$ can be analytically generated~\cite[Theorem~3.14]{Mansour-Schork-B2016} by~\eqref{Stirl-No-First-GF}.
The generating function $\bigl[\frac{\ln(1+x)}{x}\bigr]^k$ can be regarded as the power of the normalized remainder $(T_0[\ln(1+x)])^k$.
\par
In~\cite{Howard-Fib-1980i4}, Howard defined $s_r(k,\ell)$ by
\begin{equation}\label{1st-Stirl-gen-log}
\Biggl(\ln\frac{1}{1-x}-\sum_{k=1}^{r}\frac{x^k}{k}\Biggr)^\ell
=\ell!\sum_{k=(r+1)\ell}^{\infty}s_r(k,\ell)\frac{x^k}{k!}, \quad r\in\mathbb{N}_0.
\end{equation}
It is clear that $s_0(k,\ell)=(-1)^{k+\ell}s(k,\ell)$.
\par
The equation~\eqref{1st-Stirl-gen-log} can be reformulated as
\begin{multline}\label{1st-Stirl-gen-log-RF}
\Biggl[(-1)^r\frac{r+1}{x^{r}}\Biggl(\frac{\ln(1+x)}{x}-\sum_{k=0}^{r-1}(-1)^k\frac{x^{k}}{k+1}\Biggr)\Biggr]^\ell\\
=\frac{\ell!(r+1)^\ell}{[(r+1)\ell]!} \sum_{k=0}^{\infty}(-1)^k\frac{s_r(k+(r+1)\ell,\ell)} {\binom{k+(r+1)\ell}{k}}\frac{x^{k}}{k!}
\end{multline}
for $r\in\mathbb{N}_0$. The function
\begin{multline*}
(-1)^{r}\frac{r+1}{x^{r}}\Biggl[\frac{\ln(1+x)}{x}-\sum_{k=0}^{r-1}(-1)^{k}\frac{x^{k}}{k+1}\Biggr]\\
=(-1)^{r}\frac{r+1}{x^{r+1}}\Biggl[\ln(1+x)-\sum_{k=1}^{r}(-1)^{k-1}\frac{x^{k}}{k}\Biggr]
\end{multline*}
in~\eqref{1st-Stirl-gen-log-RF} is just the normalized remainder $T_{r-1}\bigr[\frac{\ln(1+x)}{x}\bigr]$ for $r\in\mathbb{N}_0$ or the normalized remainder $T_r[\ln(1+x)]$; see also~\cite[Remarks~5 and~8]{Spivey-JIMA.tex} and~\cite[Remarks~5 and~8]{Spivey-ID-AP-JIMA.tex}.

\subsection{A Generating Function of Broder's Generalization of Stirling Numbers of the First Kind}
Theorem~15 in~\cite{Broder-1984} reads that the $r$-Stirling numbers of the first kind for $r\in\mathbb{N}_0$ have the ``vertical'' exponential generating function
\begin{equation}\label{r-stirling-first}
\sum_{k}\begin{bmatrix}
k+r\\m+r
\end{bmatrix}_r\frac{x^k}{k!}=
\begin{dcases}
\frac{1}{m!}\biggl(\frac{1}{1-x}\biggr)^r\biggl(\ln\frac{1}{1-x}\biggr)^m, &m\ge0;\\
0, & m<0.
\end{dcases}
\end{equation}
We can rewrite~\eqref{r-stirling-first} in the form
\begin{equation}\label{r-stirling-first-reform}
\begin{split}
\biggl(\frac{1}{1+x}\biggr)^r\biggl[\frac{\ln(1+x)}{x}\biggr]^m
&=\biggl(\frac{1}{1+x}\biggr)^r\biggl(T_0\biggr[\frac{\ln(1+x)}{x}\biggr]\biggr)^m\\
&=\sum_{k=0}^\infty(-1)^k\frac{\left[\begin{smallmatrix}
k+m+r\\m+r
\end{smallmatrix}\right]_r}{\binom{k+m}{m}}\frac{x^k}{k!}, \quad r,m\in\mathbb{N}_0
\end{split}
\end{equation}
for $|x|<1$.
Taking $r=0$ in~\eqref{1st-Stirl-gen-log-RF} and~\eqref{r-stirling-first-reform} and comparing with~\eqref{Stirl-No-First-GF} give
\begin{equation*}
\begin{bmatrix}
k\\m
\end{bmatrix}_0=s_0(k,m)=(-1)^{k+m}s(k,m), \quad k,m\in\mathbb{N}_0.
\end{equation*}
See also~\cite[Remark~6]{Spivey-JIMA.tex} and~\cite[Remark~6]{Spivey-ID-AP-JIMA.tex}.

\subsection{A Generating Function of Broder--Carlitz's Generalization of Stirling Numbers of the Second Kind}\label{sec4.5}
Theorem~16 in~\cite{Broder-1984} and Eq.~(3.9) in~\cite{Carlitz-Fibonacci-1980-I} state that the $r$-Stirling numbers of the second kind for $r\in\mathbb{N}_0$ have the exponential generating function
\begin{equation}\label{r-stirling-second}
\sum_{k}\left\{\begin{matrix}
k+r\\ m+r
\end{matrix}\right\}_r\frac{x^k}{k!}=
\begin{dcases}
\frac{1}{m!}\te^{r x}(\te^x-1)^m, & m\ge0;\\
0, & m<0.
\end{dcases}
\end{equation}
We can reformulate~\eqref{r-stirling-second} in the form
\begin{equation}\label{small-matrix-stirl}
\te^{r x}(T_0[\te^x])^m=\sum_{k=0}^\infty\frac{\left\{\begin{smallmatrix}
k+m+r\\ m+r
\end{smallmatrix}\right\}_r}{\binom{k+m}{m}}\frac{x^k}{k!}, \quad r,m\in\mathbb{N}_0, \quad |x|<\infty.
\end{equation}
Taking $r=0$ in~\eqref{T-Associate-Stirling} and~\eqref{small-matrix-stirl} and comparing with~\eqref{2Stirl-funct-rew} yield
\begin{equation*}
\left\{\begin{matrix}
k\\ m
\end{matrix}\right\}_0
=S_0(k,m)=S(k,m), \quad k,m\in\mathbb{N}_0.
\end{equation*}
See also~\cite[Remark~7]{Spivey-JIMA.tex} and~\cite[Remark~7]{Spivey-ID-AP-JIMA.tex}.

\subsection{Unified Treatments of Generating Functions of Bernoulli Numbers and Stirling Numbers}
In the paper~\cite{Gottfried-Qi.tex}, by virtue of a determinantal formula for derivatives of the ratio between two differentiable functions, in view of the Fa\`a di Bruno formula, and with help of several identities and closed-form formulas for partial Bell polynomials $\bell_{n,k}$ defined by~\eqref{Bell2nd-Dfn-Eq}, we
\begin{enumerate}[1)]
\item
established thirteen Maclaurin series expansions of the functions
\begin{gather*}
\ln\frac{\te^x+1}{2}, \quad \ln\frac{\te^x-1}{x}, \quad \ln\cosh x, \\
\ln\frac{\sinh x}{x}, \quad \biggl[\frac{\ln(1+x)}{x}\biggr]^r, \quad \biggl(\frac{\te^x-1}{x}\biggr)^r
\end{gather*}
for $r=\pm\frac{1}{2}$ or $r\in\mathbb{R}$ in terms of the Dirichlet eta function $\eta(1-2k)$, the Riemann zeta function $\zeta(1-2k)$, and the Stirling numbers of the first and second kinds $s(n,k)$ and $S(n,k)$, respectively;
\item
presented four determinantal expressions and three recursive relations for the Bernoulli numbers $B_{2n}$;
\item
found out three closed-form formulas for the Bernoulli numbers $B_{2n}$ and the generalized Bernoulli numbers $B_n^{(r)}$ in terms of the Stirling numbers of the second kind $S(n,k)$;
\item
deduced two combinatorial identities for the Stirling numbers of the second kind $S(n,k)$;
\item
acquired two combinatorial identities, which can be regarded as diagonal recursive relations, involving the Stirling numbers of the first and second kinds $s(n,k)$ and $S(n,k)$;
\item
recovered an integral representation and a closed-form formula, and establish an alternative explicit and closed-form formula, for the Bernoulli numbers of the second kind $b_n$ in terms of the Stirling numbers of the first kind $s(n,k)$;
\item
obtained three identities connecting the Stirling numbers of the first and second kinds $s(n,k)$ and $S(n,k)$.
\end{enumerate}
The highlights of the paper~\cite{Gottfried-Qi.tex} include the unification $\bigl(\frac{\te^x-1}{x}\bigr)^r$ of the generating functions for the Bernoulli numbers $B_n$ and the Stirling numbers of the second kind $S(n,k)$, the unification $\bigl[\frac{\ln(1+x)}{x}\bigr]^r$ of the generating functions for the Bernoulli numbers of the second kind $b_n$, which can be generated~\cite{Davis-1957-Monthly, Bernoulli-Stirling-Beograd.tex, Bernoulli2nd-Property.tex} by
\begin{equation*}
\frac{x}{\ln(1+x)}=\sum_{n=0}^\infty b_nx^n, \quad |x|<1,
\end{equation*}
and the Stirling numbers of the first kind $s(n,k)$, and the disclosure of the transformations between these two unifications.

\subsection{Revisiting Generating Function of Genocchi Numbers}
The Genocchi numbers $G_n$ are generated~\cite[p.~49]{Comtet-Combinatorics-74} by
\begin{equation*}
\frac{2x}{\te^x+1}=\sum_{n=1}^{\infty}G_n\frac{x^n}{n!},
\end{equation*}
which can be reformulated as
\begin{equation*}
\frac{1}{T_{-1}[\te^x+1]}=\sum_{n=0}^{\infty}G_{n+1}\frac{x^{n}}{(n+1)!}.
\end{equation*}

\subsection{Revisiting Generating Functions of Euler Numbers and Polynomials}
The Euler numbers and polynomials $E_n$ and $E_n(t)$ are defined~\cite[p.~48]{Comtet-Combinatorics-74} respectively by
\begin{equation*}
\frac{2\te^x}{\te^{2x}+1}=\sum_{k=0}^{\infty}E_n\frac{x^n}{n!}, \quad |x|<\frac{\pi}{2} \quad\text{and}\quad
\frac{2\te^{t x}}{\te^x+1}=\sum_{k=0}^{\infty}E_n(t)\frac{x^n}{n!}, \quad |x|<\pi.
\end{equation*}
Since
\begin{equation*}
\te^x+1=2+\sum_{n=1}^{\infty}\frac{x^n}{n!}, \quad x\in\mathbb{R},
\end{equation*}
we can regard the generating functions $\frac{2\te^x}{\te^{2x}+1}$ and $\frac{2\te^{t x}}{\te^x+1}$ as
\begin{equation*}
\frac{\te^x}{T_{-1}[\te^{2x}+1]} \quad\text{and}\quad \frac{\te^{t x}}{T_{-1}[\te^x+1]}
\end{equation*}
respectively. Therefore, we can generalize the Euler polynomials $E_n(t)$ by
\begin{equation}\label{Enksigma)}
\frac{\te^{t x}}{(T_n[\te^x+1])^\sigma}=\frac{\te^{t x}}{(T_n[\te^x])^\sigma}=\sum_{k=0}^{\infty}E_{n,k}^{(\sigma)}(t)\frac{x^k}{k!}, \quad n\in\mathbb{N}_0.
\end{equation}

\section{Basic Properties of Normalized Remainders}\label{sec-basic-property}

In this section, we recite some basic properties of normalized remainders.

\begin{thm}[{\cite[Theorem~1]{arcsine-liu-qi-mdpi.tex}}]\label{Theorem1arcsine-liu-qi-mdpi.tex}
If $T_n[f(x),x_0]$ for some $n\in\mathbb{N}_0$ exists, then the normalized remainder $T_n[f(x),x_0]$ satisfies
\begin{equation*}
T_n[\alpha+f(x),x_0]=T_n[f(x),x_0], \quad n\in\mathbb{N}_0,
\end{equation*}
where $\alpha$ is a real number.
\end{thm}

\begin{thm}[{\cite[Theorem~1]{LiuXL-Arc-three.tex}}]\label{One-more-simple-property}
If $T_n[f(x),x_0]$ for some $n\in\mathbb{N}_0$ exists, then
\begin{equation*}
T_n[\beta f(x),x_0]=T_n[f(x),x_0], \quad n\in\mathbb{N}_0,
\end{equation*}
where $\beta\ne0$ is a real number.
\end{thm}

\begin{thm}[{\cite[Theorem~1]{arcsine-liu-qi-mdpi.tex}}]\label{nR-power}
If $f^{(i)}(x_0)=0$ for $0\le i<k$ and the normalized remainder $T_n[f(x),x_0]$ exists for some $n\ge k$, then the normalized remainder $T_n[f(x),x_0]$ satisfies
\begin{equation}\label{nR-power-eq}
T_n[f(x),x_0]=T_{n-k}\biggl[\frac{f(x)}{(x-x_0)^k},x_0\biggr], \quad n\ge k\in\mathbb{N}_0.
\end{equation}
\end{thm}

Utilizing Theorems~\ref{Theorem1arcsine-liu-qi-mdpi.tex}, \ref{One-more-simple-property}, and~\eqref{nR-power} with the relations
\begin{equation*}
\arccos z=\frac{\pi}{2}-\arcsin z, \quad \arctan z=\frac{\pi}{2}-\arccot z,
\end{equation*}
we arrive at
\begin{align*}
T_{2n}[\arcsin x]&=T_{2n-1}\biggl[\frac{\arcsin x}{x}\biggr], & T_{2n}[\arctan x]&=T_{2n-1}\biggl[\frac{\arctan x}{x}\biggr],\\
T_{2n}[\arccos x]&=T_{2n}[\arcsin x], & T_{2n}[\arctan x]&=T_{2n}[\arccot x],
\end{align*}
and
\begin{align*}
T_{2n}\biggl[\arccos x-\frac{\pi}{2}\biggr]&=T_{2n-1}\biggl[\frac{\arccos x-\frac{\pi}{2}}{x}\biggr],\\
T_{2n}\biggl[\arccot x-\frac{\pi}{2}\biggr]&=T_{2n-1}\biggl[\frac{\arccot x-\frac{\pi}{2}}{x}\biggr].
\end{align*}
See also the paper~\cite{LiuXL-Arc-three.tex}.
\par
The relations~\eqref{ln-sec-cos-tail} and~\eqref{tan-sec-square} can also be deduced from Theorems~\ref{Theorem1arcsine-liu-qi-mdpi.tex} and~\ref{One-more-simple-property}, respectively.

\begin{thm}
If $f(x)$ is an absolutely monotonic function on an interval $[0,a]\subset[0,\infty)$ and its normalized remainder $T_n[f(x)]$ exists for some $n\in\mathbb{N}_0$, then its normalized remainder $T_n[f(x)]$ for some $n\in\mathbb{N}_0$ is also an absolutely monotonic function on $[0,a]\subset[0,\infty)$.
\end{thm}

\begin{proof}
This follows from combining Definition~\ref{defn-normaliz} with a statement in~\cite[p.~156]{Widder-Laplace-Transform-41} which reads that, if $f(x)$ is absolutely monotonic in $[0,a]$ and if $f(x)=\sum_{n=0}^{\infty}A_nx^n$ for $x\in[0,a]$, then all coefficients $A_n$ are non-negative.
\end{proof}

\section{Abu-Ghuwaleh's work}
In the preprint~\cite{Mohammad-Abu-Ghuwaleh}, Abu-Ghuwaleh advanced the abstract and general theory of the normalized remainders, developing it from a fresh angle and perspective. It offers an intrinsic and innovative treatment of the topic and opens up an exciting new direction of research.
\par
In the preprint~\cite{Mohammad-Abu-Ghuwaleh}, Abu-Ghuwaleh worked with analytic functions whose Maclaurin coefficients never vanish and gave a variant of Definition~\ref{defn-normaliz}.

\begin{defn}[{\cite[Definition~1]{Mohammad-Abu-Ghuwaleh}}]\label{Definition1-Mohammad}
For $R\in(0,\infty]$, let $\mathscr{A}_R^\times$ denote the set of all analytic functions
\begin{equation*}
f(z)=\sum_{n=0}^{\infty}a_nz^n\in\mathcal{O}(D_R)
\end{equation*}
such that $a_n\ne0$ for all $n\ge0$, where
\begin{equation*}
D_R=\{z\in\mathbb{C}:|z|<R\}
\end{equation*}
is the open disc on the complex plane $\mathbb{C}$ with radius $R$ and $\mathcal{O}(D_R)$ denotes the set of all analytic functions defined on the disc $D_R$.
\par
For $f\in\mathscr{A}_R^\times$, define the normalized Maclaurin tail of order $n$ by
\begin{equation*}
T_n^f(z)=\sum_{k=0}^{\infty}\frac{a_{n+k}}{a_n}z^k, \quad n\in\mathbb{N}_0.
\end{equation*}
\end{defn}

\begin{rem}
It is easy to see that
$$
T_0^f(z)=\frac{f(z)}{a_0}=\frac{f(z)}{f(0)}.
$$
\par
For $f\in\mathscr{A}_R^\times$, it is obvious that
\begin{equation*}
a_n=\frac{f^{(n)}(0)}{n!}, \quad n\in\mathbb{N}_0.
\end{equation*}
Replacing $\frac{f^{(n)}(0)}{n!}$ by $a_n$ and $z_0$ by $0$ in Definition~\ref{defn-normaliz}, the equation~\eqref{T(f(x))} becomes
\begin{align*}
T_n[f(z)]&=
\begin{dcases}
\frac1{a_{n+1}z^{n+1}} \Biggl(\sum_{k=0}^{\infty}a_kx^k-\sum_{k=0}^{n} a_kx^k\Biggr), & z\in (-R,R)\setminus\{0\}\\
1, & z=0
\end{dcases}\\
&=
\begin{dcases}
\frac1{a_{n+1}z^{n+1}}\sum_{k=n+1}^{\infty}a_kx^k, & z\in (-R,R)\setminus\{0\}\\
1, & z=0
\end{dcases}\\
&=
\begin{dcases}
\sum_{k=0}^{\infty}\frac{a_{n+k+1}}{a_{n+1}}z^{k}, & z\in (-R,R)\setminus\{0\}\\
1, & z=0
\end{dcases}
\end{align*}
for $n\in\mathbb{N}_0$. This implies that
\begin{equation}\label{TfTn}
T_{n}^f(z)=T_{n-1}[f(z)], \quad n\in\mathbb{N}_0.
\end{equation}
Consequently, Definition~\ref{Definition1-Mohammad} is essentially the same as Definition~\ref{defn-normaliz}.
\end{rem}

\begin{rem}
The requirement $a_n\ne0$ for $n\in\mathbb{N}_0$ in Definition~\ref{defn-normaliz} excludes all odd functions from objects of study. For an even function $f(z^2)$, by a transform $z^2\to w$, the function $f(w)$ is included in objects of study.
\end{rem}

We now recite the main results in~\cite{Mohammad-Abu-Ghuwaleh} and give several remarks on them.

\begin{thm}[{Radius invariance~\cite[Proposition~1]{Mohammad-Abu-Ghuwaleh}}]
Let $f\in\mathscr{A}_R^\times$. Then the power series $T_n^f$ has the same radius of convergence. In particular, $T_n^f\in\mathscr{G}_R^\times$, where $\mathscr{G}_R^\times$ denotes the set of all analytic functions
\begin{equation}\label{F(z)}
F(z)=1+\sum_{n=1}^{\infty}c_nz^n\in\mathcal{O}(D_R)
\end{equation}
such that $c_n\ne0$ for $n\in\mathbb{N}$.
\end{thm}

\begin{rem}
Since
\begin{equation*}
T_n^f(z)=\frac{1}{a_n}\sum_{k=0}^{\infty}a_{n+k}z^k, \quad n\in\mathbb{N}
\end{equation*}
the radius invariance property is almost obvious.
\par
Since $T_n^f(0)=1$ and $a_n\ne0$ for $n\in\mathbb{N}_0$, the relation $T_n^f\in\mathscr{G}_R^\times$ for $n\in\mathbb{N}_0$ is trivial.
\end{rem}

\begin{rem}
The convergent radius of the power series expansion~\eqref{Bernoulli-GF} is $2\pi$.
However, we established in the paper~\cite{Bernouli-No-Tail.tex} some properties of the normalized remainder $T_{2n+1}\bigl[\frac{x}{\te^x-1}\bigr]$ for $n\in\mathbb{N}_0$ on the real line $\mathbb{R}$. For details, please refer to Section~\ref{sec3.31} in this chapter.
\end{rem}

\begin{thm}[{Normalized left shift~\cite[Proposition~2]{Mohammad-Abu-Ghuwaleh}}]
If $F\in\mathscr{G}_R^\times$, then
\begin{equation*}
\mathcal{S}(F)(z)=\sum_{k=0}^{\infty}\frac{c_{k+1}}{c_1}z^k,
\end{equation*}
where the normalization operator $\mathcal{S}$ is defined for $R\in(0,\infty]$ by
\begin{equation*}
\mathcal{S}:\mathscr{G}_R^\times\to\mathscr{G}_R^\times, \qquad \mathcal{S}(F)(z)=\frac{F(z)-1}{zF'(0)}.
\end{equation*}
Consequently, we have
\begin{equation}\label{SnF(z)}
\mathcal{S}^n(F)(z)=\sum_{k=0}^{\infty}\frac{c_{n+k}}{c_n}z^k, \quad n\in\mathbb{N}.
\end{equation}
Moreover, the operator $\mathcal{S}$ is surjective but not injective. More precisely, if $F\in\mathscr{G}_R^\times$, then the full fiber over $F$ is
\begin{equation*}
\bigl\{G\in\mathscr{G}^\times:\mathcal{S}(G)=F\bigr\}=\bigl\{1+\lambda wF(w):\lambda\in\mathbb{C}^\times=\mathbb{C}\setminus\{0\}\bigr\}.
\end{equation*}
\end{thm}

\begin{rem}
It is clear that, if $f\in\mathscr{A}_R^\times$, then $T_{n-1}[f(z)]\in\mathscr{G}_R^\times$ for all $n\in\mathbb{N}_0$.
\par
It is not difficult to verify that
\begin{equation}\label{S-T}
\mathcal{S}(F)(z)=T_0[F(z)].
\end{equation}
Hence, we have
\begin{equation*}
\mathcal{S}^2(F)(z)=\mathcal{S}(T_0[F(z)])(z)=T_0[T_0[F(z)]]=T_0^2[F(z)].
\end{equation*}
Inductively, we obtain
\begin{equation*}
\mathcal{S}^n(F)(z)=\underbrace{T_0[T_0[\dotsm T_0[T_0[T_0}_n[F(z)]]]]]=T_0^n[F(z)], \quad n\in\mathbb{N}.
\end{equation*}
\end{rem}

\begin{thm}[{Exact normalization law~\cite[Theorem~1]{Mohammad-Abu-Ghuwaleh}}]
Let $f\in\mathscr{A}_R^\times$ and $n\in\mathbb{N}_0$. Then
\begin{equation}\label{T(n+1)T(n)}
T_{n+1}^f=\mathcal{S}\bigl(T_n^f\bigr).
\end{equation}
Equivalently,
\begin{equation*}
T_n^f=\mathcal{S}^n\bigl(T_0^f\bigr).
\end{equation*}
In coefficient form,
\begin{equation}\label{T(n+1)T(n)-3}
T_n^f(z)=1+\frac{a_{n+1}}{a_n}z T_{n+1}^f(z).
\end{equation}
\end{thm}

\begin{rem}
By the relations~\eqref{TfTn} and~\eqref{S-T}, we can rewrite~\eqref{T(n+1)T(n)} to~\eqref{T(n+1)T(n)-3} as
\begin{align*}
T_{n}[f(z)]&=T_0[T_{n-1}[f(z)]],\\
T_{n-1}[f(z)]&=\underbrace{T_0[T_0[\dotsm T_0[T_0[T_0}_n[T_{-1}[f(z)]]]]]]\\
&=T_0^n[T_{-1}[f(z)]],\\
\intertext{and}
T_{n-1}[f(z)]&=1+\frac{a_{n+1}}{a_n}z T_{n}[f(z)]
\end{align*}
for $n\in\mathbb{N}_0$.
\end{rem}

Due to space limitations, we no longer review Abu‑Ghuwaleh's work in the preprints~\cite{Mohammad-Abu-Ghuwaleh} and tens of subsequent papers.

\section{Conjectures and Problems on Normalized Remainders}
In this section, we collect several conjectures and problems posed while investigating normalized remainders of several elementary functions.

\subsection{Conjectures on Normalized Remainders}

The following are nine conjectures on concrete normalized remainders.

\subsubsection{First Conjecture}
In~\cite[Remark~4.1]{ratio-NR-squ-tan.tex}, we proposed the conjecture: The ratio $\frac{T_{2n+1}[\tan^2x]}{T_{2n-1}[\tan^2x]}$ for $n\in\mathbb{N}$ is concave in $x\in\bigl(-\frac{\pi}{2},\frac{\pi}{2}\bigr)$.

\subsubsection{Second Conjecture}
In the paper~\cite[Section~4.1]{JMI-5203.tex}, we posed the following conjectures:
For $n\in\mathbb{N}$, the ratio $\frac{T_{2n-1}[\tan^2x]}{T_{2n+1}[\tan^2x]}$ is completely monotonic in $x\in\bigl(-\frac{\pi}{2},0\bigr)$ and absolutely monotonic in $x\in\bigl(0,\frac{\pi}{2}\bigr)$.

\subsubsection{Third Conjecture}
In the paper~\cite[Section~4.2]{JMI-5203.tex}, we posed the following conjectures:
For $n\in\mathbb{N}_0$, the normalized remainder $T_n[\te^x]$ is a logarithmically absolutely monotonic function in $x\in\mathbb{R}$.

\subsubsection{Fourth Conjecture}
In the paper~\cite[Section~4.2]{JMI-5203.tex}, we posed the following conjectures:
For $n\in\mathbb{N}_0$, the ratio $\frac{T_{n}[\te^x]}{T_{n+1}[\te^x]}$ is an absolutely monotonic function in $x\in\mathbb{R}$.

\subsubsection{Fifth Conjecture}
In~\cite[Remark~7]{exp-norm-tail-ratio.tex}, we conjectured that the functions
\begin{equation}\label{Gn(u)3series-dfn}
G_n(x)=\sum_{j=0}^{\infty}\frac{1}{\binom{j+n+1}{n+1}}\frac{x^{j}}{j!} \sum_{j=0}^{\infty}\frac{1}{\binom{j+n+1}{n}}\frac{x^{j}}{j!} -\sum_{j=0}^{\infty}\frac{1}{\binom{j+n+2}{n+1}}\frac{x^{j}}{j!} \sum_{j=0}^{\infty}\frac{1}{\binom{j+n}{n}}\frac{x^{j}}{j!}
\end{equation}
and
\begin{equation}\label{Hn(u)3series-dfn}
H_n(x)=\sum_{j=0}^{\infty}\frac{1}{\binom{j+n+2}{n}}\frac{x^{j}}{j!} \sum_{j=0}^{\infty}\frac{1}{\binom{j+n}{n}}\frac{x^{j}}{j!} -\Biggl[\sum_{j=0}^{\infty}\frac{1}{\binom{j+n+1}{n}}\frac{x^{j}}{j!}\Biggr]^2
\end{equation}
for $n\in\mathbb{N}$ are absolutely monotonic in $x\in\mathbb{R}$.
\par
The absolute monotonicity of $G_n(x)$ and $H_n(x)$ are stronger than the inequalities~\eqref{log-exp-last-proof-ineq} and~\eqref{log-exp-last-inverse-ineq}.
\par
In~\cite[Section~1]{ser-int-norm-tail.tex}, we conjectured that the function
\begin{equation}\label{Jn(u)3series-dfn}
J_n(x)=\Biggl[\sum_{j=0}^{\infty}\frac{1}{\binom{j+n+2}{n+1}}\frac{x^{j}}{j!}\Biggr]^2 -\frac{n+1}{n+2}\sum_{j=0}^{\infty}\frac{1}{\binom{j+n+3}{n+2}}\frac{x^{j}}{j!} \sum_{j=0}^{\infty}\frac{1}{\binom{j+n+1}{n}}\frac{x^{j}}{j!}
\end{equation}
for $n\in\mathbb{N}_0$ is absolutely monotonic in $x\in\mathbb{R}$,

\subsubsection{Sixth Conjecture}
Since the quantity $\mathcal{Q}(q,2\ell)$ is positive for $q,\ell\in\mathbb{N}$, we conjectured in~\cite[Section~1]{arcsine-liu-qi-mdpi.tex} that for $q\in\mathbb{N}$ and $n\in\mathbb{N}_0$ the normalized remainder $T_{2n+1}\bigl[\bigl(\frac{\arcsin x}{x}\bigr)^q\bigr]$ is logarithmically absolutely monotonic in $x\in[0,1)$.

\subsubsection{Seventh Conjecture}
In the paper~\cite{LiuXL-Arc-three.tex}, we conjectured that
the normalized remainder $T_{2n}[\arctan x]$ for $n\in\mathbb{N}_0$ is logarithmically concave on $(0,1)$.

\subsubsection{Eighth Conjecture}
In the paper~\cite{LiuXL-Arc-three.tex}, we conjectured that the ratio $\frac{T_{2n}[\arcsin x]}{T_{2n+2}[\arcsin x]}$ for $n\in\mathbb{N}_0$ is concave in $x\in(-1,1)$.

\subsubsection{Ninth Conjecture}
In the paper~\cite{LiuXL-Arc-three.tex}, we conjectured that
the ratio $\frac{T_{2n}[\arctan x]}{T_{2n+2}[\arctan x]}$ for $n\in\mathbb{N}_0$ is increasing in $x\in(0,\infty)$.

\subsection{Problems on Normalized Remainders}
The following are eight problems related to normalized remainders.

\subsubsection{First Problem}
Expand the logarithm $\ln T_{2n}[\tan x]$ into a Maclaurin series.

\subsubsection{Second Problem}
Let $L_n(x)=\ln T_{2n-2}[\tan x]$ for $n\in\mathbb{N}$ and let
\begin{equation*}
Q_{m,n}(x)=
\begin{dcases}
\frac{L_n(x)}{L_m(x)},\quad 0<|x|<\frac{\pi}{2}\\
\frac{(m+1)(2m+1)\bigl(2^{2m}-1\bigr)(2^{2n+2}-1)B_{2m}B_{2n+2}} {(n+1)(2n+1)(2^{2m+2}-1)(2^{2n}-1)B_{2m+2}B_{2n}},&x=0
\end{dcases}
\end{equation*}
for $n>m\in\mathbb{N}$. It is clear that $Q_{1,2}(x)=\frac{\ln F(x)}{\ln G(x)}$, where $F(x)$ and $G(x)$ are defined in~\eqref{f(x)-def-cases} and~\eqref{G(x)Tangent}.
\par
Show that the function $Q_{m,n}(x)$ for $m>n\in\mathbb{N}$ is decreasing on $\bigl(-\frac{\pi}{2},\frac{\pi}{2}\bigr)$.

\subsubsection{Third Problem}
Expand the logarithm $\ln T_{2n}[\sin x]$ for $n\ge2$ into a Maclaurin series.

\subsubsection{Fourth Problem}
Discuss the monotonicity and convexity of the ratio $\frac{\ln T_{2n+1}[\cos x]}{\ln T_{2n-1}[\cos x]}$ for $n\in\mathbb{N}_0$ on $\bigl[0,\frac{\pi}2\bigr]$.

\subsubsection{Fifth Problem}
In~\cite[Remark~6]{Spivey-JIMA.tex} and~\cite[Remark~8]{Spivey-ID-AP-JIMA.tex}, we proposed the problem: Investigate properties of the sequence $F(r,s,m,k)$ generated by
\begin{align*}
\biggl(\frac{1}{1+x}\biggr)^r(T_s[\ln(1+x)])^m
&=\biggl(\frac{1}{1+x}\biggr)^r\biggl(T_s\biggr[\frac{\ln(1+x)}{x}\biggr]\biggr)^m\\
&=\sum_{k=0}^\infty F(r,s,m,k)\frac{x^k}{k!}
\end{align*}
for $r,m\in\mathbb{C}$, $s\in\mathbb{N}_0$, and $|x|<1$.
\par
The integral representation~\eqref{Tnexp(x)-int} and the ideas with techniques in~\cite[Section~7]{Gottfried-Qi.tex} should be useful when studying this problem.

\subsubsection{Sixth Problem}
In~\cite[Remark~7]{Spivey-JIMA.tex} and~\cite[Remark~7]{Spivey-ID-AP-JIMA.tex}, we proposed the problem: Generally, we can consider a new kind of polynomials $A_{n,k}^{(\sigma)}(t)$ generated by
\begin{equation*}
\frac{\te^{x t}}{(T_{n-1}[\te^x])^\sigma}
=\sum_{k=0}^{\infty}A_{n,k}^{(\sigma)}(t)\frac{x^k}{k!}, \quad n\in\mathbb{N}, \quad \sigma,x,t\in\mathbb{C}.
\end{equation*}
It is obvious that the polynomial sequence $A_{n,k}^{(\sigma)}(t)$ is a unified generalization of $B_k$, $B_k(t)$, $A_k$, $A_k(t)$, $A_{n,k}(t)$, $A_{n,k}(t)$, and the Stirling numbers of the second kind $S(n,k)$, and the Stirling polynomials $S_k(x)$.
\par
It is not difficult to see that the sequences $A_{n,k}^{(\sigma)}(t)$ and $E_{n,k}^{(\sigma)}(t)$, generated in~\eqref{Enksigma)}, are essentially the same one.
\par
The integral representation~\eqref{Tnexp(x)-int} should be useful when studying this problem.

\subsubsection{Seventh Problem}
It is well-known that, if a function $f(x)$ is of the period $2\ell>0$ on $\mathbb{R}$ and the integrals
\begin{equation*}
a_n=\frac{1}{\ell}\int_{-\ell}^{\ell}f(x)\cos\frac{n\pi x}{\ell}\td x, \quad n\in\mathbb{N}_0
\end{equation*}
and
\begin{equation*}
b_n=\frac{1}{\ell}\int_{-\ell}^{\ell}f(x)\sin\frac{n\pi x}{\ell}\td x, \quad n\in\mathbb{N}
\end{equation*}
exist, then the trigonometric series
\begin{equation*}
\frac{a_0}{2}+\sum_{n=1}^{\infty}\biggl(a_n\cos\frac{n\pi x}{\ell}+b_n\sin\frac{n\pi x}{\ell}\biggr)
\end{equation*}
is called the Fourier series of the function $f(x)$.
\par
Similar to Definition~\ref{defn-normaliz}, consider the function
\begin{equation*}
F_n[\ell,f(x)]=\frac{\frac{a_0}{2}+\sum_{k=1}^{n}\bigl(a_k\cos\frac{k\pi x}{\ell}+b_k\sin\frac{k\pi x}{\ell}\bigr)} {a_{n+1}\cos\frac{(n+1)\pi x}{\ell}+b_{n+1}\sin\frac{(n+1)\pi x}{\ell}}
\end{equation*}
for $n\in\mathbb{N}$ and $x\in\mathbb{R}$ and discuss its properties.

\subsubsection{Eighth Problem}
Let $\boldsymbol{s}_n=(s_1,s_2,\dotsc, s_{n})$ and $\boldsymbol{t}_n=(t_1,t_2,\dotsc, t_{n})\in\mathbb{R}^{n}$.
A $n$-tuple $\boldsymbol{s}_n$ is said to strictly majorize $\boldsymbol{t}_n$, denoted by $\boldsymbol{t}_n\prec\boldsymbol{s}_n$ or $\boldsymbol{s}_n\succ\boldsymbol{t}_n$, if
\begin{equation*}
\bigl(s_{[1]},s_{[2]},\dotsc, s_{[n]}\bigr)\ne\bigl(t_{[1]},t_{[2]},\dotsc, t_{[n]}\bigr),\qquad
\sum_{i=1}^n s_i=\sum_{i=1}^n t_i,
\end{equation*}
and
\begin{equation*}
\sum_{i=1}^k s_{[i]}>\sum_{i=1}^k t_{[i]}, \quad 1\le k\le n-1
\end{equation*}
are valid, where
$$
s_{[1]}\ge s_{[2]}\ge\dotsm \ge s_{[n]}
\quad \text{and}\quad
t_{[1]}\ge t_{[2]}\ge\dotsm \ge t_{[n]}
$$
are rearrangements of $\boldsymbol{s}_n$ and $\boldsymbol{t}_n$ in a descending order; see~\cite[p.~8, Defnition~A.1]{Marshall-Olkin-Arnold}.
\par
In~\cite[Section~5]{ser-int-norm-tail.tex}, the function $G_n(x)$ defined by~\eqref{Gn(u)3series-dfn} was generalized as
\begin{equation}\label{Gmnk-ell(u)}
G_{m,n;k,\ell}(x)=-\begin{vmatrix}
S_{m,n}(x) & S_{m,n+\ell}(x)\\
S_{m+k,n}(x) & S_{m+k,n+\ell}(x)
\end{vmatrix}
\end{equation}
for $m,n\in\mathbb{N}_0$ and $k,\ell\in\mathbb{N}$, where $S_{m,n}(x)$ is defined by~\eqref{S(mn)(u)-ser-defn}.
\par
In~\cite[Section~5]{ser-int-norm-tail.tex}, basing on an inequality in~\cite[Theorem~7]{ser-int-norm-tail.tex}, we defined
\begin{equation}\label{K-ell(u)}
K_n(x)=\frac{n+1}{n+2}S_{1,n+2}(x) S_{1,n}(x) -\frac{n}{n+1}S_{1,n+1}^2(x)
\end{equation}
for $n\in\mathbb{N}_0$ and $x\in\mathbb{R}$.
\par
The functions $H_n(x)$ defined by~\eqref{Hn(u)3series-dfn}, $J_n(x)$ defined by~\eqref{Jn(u)3series-dfn}, and $K_n(x)$ defined by~\eqref{K-ell(u)} can be generalized as follows:
\begin{enumerate}[1)]
\item
For $n\in\mathbb{N}_0$ and $(s_1,s_2),(t_1,t_2)\in\mathbb{N}_0^2$ such that $(s_1,s_2)\succ(t_1,t_2)$, define
\begin{equation}\label{Hnboldsymbols2boldsymbolt2(u)}
H_{n;\boldsymbol{s}_2,\boldsymbol{t}_2}(x)=S_{s_1,n+1}(x) S_{s_2,n+1}(x) -S_{t_1,n+1}(x)S_{t_2,n+1}(x).
\end{equation}
\item
For $m\in\mathbb{N}_0$ and $(s_1,s_2),(t_1,t_2)\in\mathbb{N}_0^2$ such that $(s_1,s_2)\succ(t_1,t_2)$, define
\begin{equation}\label{Jmboldsymbols2boldsymbolt2(u)}
J_{m;\boldsymbol{s}_2,\boldsymbol{t}_2}(x)=S_{m,t_1}(x)S_{m,t_2}(x) -A(m;s_1,s_2;t_1,t_2)S_{m,s_1}(x)S_{m,s_2}(x).
\end{equation}
\item
For $m\in\mathbb{N}_0$ and $(s_1,s_2),(t_1,t_2)\in\mathbb{N}_0^2$ such that $(s_1,s_2)\succ(t_1,t_2)$, define
\begin{equation}\label{Kmboldsymbols2boldsymbolt2(u)}
K_{m;\boldsymbol{s}_2,\boldsymbol{t}_2}(x)=S_{m,s_1}(x)S_{m,s_2}(x)-B(m;s_1,s_2;t_1,t_2) S_{m,t_1}(x)S_{m,t_2}(x).
\end{equation}
\end{enumerate}
\par
Problem~3 in~\cite{ser-int-norm-tail.tex} stated that,
\begin{enumerate}[1)]
\item
for $m,n\in\mathbb{N}_0$ and $k,\ell\in\mathbb{N}$, the function $G_{m,n;k,\ell}(x)$ defined by~\eqref{Gmnk-ell(u)} is absolutely monotonic in $x\in\mathbb{R}$;
\item
for $n\in\mathbb{N}_0$ and $(s_1,s_2),(t_1,t_2)\in\mathbb{N}_0^2$ such that $(s_1,s_2)\succ(t_1,t_2)$, the function $H_{n;\boldsymbol{s}_2,\boldsymbol{t}_2}(x)$ defined by~\eqref{Hnboldsymbols2boldsymbolt2(u)} is absolutely monotonic in $x\in\mathbb{R}$;
\item
for $m\in\mathbb{N}_0$ and $(s_1,s_2),(t_1,t_2)\in\mathbb{N}_0^2$ such that $(s_1,s_2)\succ(t_1,t_2)$, find the best constants $0<A(m;s_1,s_2;t_1,t_2)<1$ and $0<B(m;s_1,s_2;t_1,t_2)<1$ such that the functions $J_{m;\boldsymbol{s}_2,\boldsymbol{t}_2}(x)$ defined by~\eqref{Jmboldsymbols2boldsymbolt2(u)} and $K_{m;\boldsymbol{s}_2,\boldsymbol{t}_2}(x)$ defined by~\eqref{Kmboldsymbols2boldsymbolt2(u)} are absolutely monotonic in $x\in\mathbb{R}$.
\end{enumerate}

\section{Conclusions}

Since 2023, through the study of numerous concrete examples, Qi and his collaborators have identified a recurring pattern. Building on this observation, they introduced the concept of the normalized remainder, deliberately choosing this term and subsequently exploring its historical background and mathematical significance. They investigated its fundamental properties and sought to uncover potential applications. The formulation of such a concept demands experience, insight, courage, and intellectual acuity.
\par
In this chapter, Qi presented a synthesis of the research process and the principal findings related to the normalized remainder.
\par
In the preprint~\cite{Mohammad-Abu-Ghuwaleh}, Abu-Ghuwaleh advanced the subject in a different direction and at a more structural level by examining the broader dynamical and theoretical framework surrounding the normalized remainder family. This work further develops and formalizes the concept, firmly establishing its place within the field—precisely the goal Qi has been striving toward in his research over the past several years. Thanks to~\cite{Mohammad-Abu-Ghuwaleh}, the notion of the normalized remainder now carries a deeper and more substantial meaning.
\par
Normalized remainders constitute a new concept in mathematical analysis, combinatorial number theory, and other areas. Although novel, the idea is rooted in well-established mathematical principles and supported by a rich historical context. The main results reviewed in this chapter represent initial steps and a foundational platform for further investigation. Collectively, they demonstrate that the study of normalized remainders is both compelling and meaningful, with promising avenues for future development.

\begin{rem}
This chapter is a slightly extended, modified, and corrected version of the electronic arXiv preprints~\cite{NR-arXiv-v3, Qi-CAAO.tex, GM-VG-QiF-Ch.tex}.
\end{rem}

\paragraph{Acknowledgements}
The author is grateful to anonymous referees for their careful comments and helpful suggestions to the original version of this chapter.

\end{document}